%% file: paper.tex
\documentclass[10pt]{amsart}
\usepackage{hyperref,amsmath,amsthm,amsfonts,amscd,flafter,epsf,amssymb}
\usepackage{epsfig}
\usepackage{amsfonts}
\usepackage{amssymb}
\usepackage{amsmath}
\usepackage{amsthm}
\usepackage{latexsym}
\usepackage{amscd}
\usepackage{flafter}
\usepackage{epsf}

 \input{command1}

\input{standardcommand}

\begin{document}

\title{Heegaard Floer homology and Morse surgery}%
\author{Eaman Eftekhary}%
\address{Math Department, Harvard University, 1 Oxford Street, Cambridge, MA 02138}%
\email{eaman@math.harvard.edu}%

\thanks{The authors is partially supported by a NSF grant}%

\maketitle
\begin{abstract}
We establish surgery formulas for the filtration of the Heegaard
Floer complex associated with $\frac{p}{q}$-surgery
$Y_{\frac{p}{q}}(K)$ on a null-homologous knot $(Y,K)$, induced by
the core of the attached solid torus (which produces the surgery).
This would generalize the result of Ozsv\'ath and Szab\'o in
\cite{OS-Qsurgery}. We will also re-prove that surgery on
non-trivial knots can not produce $S^3$, as a corollary of
non-vanishing results for $\widehat{\text{HFK}}(K_{\frac{p}{q}})$
where $K$ is a knot in $S^3$.
\end{abstract}
\section{Introduction}
Suppose that $K$ is a null-homologous knot in the three-manifold
$Y$ and suppose that $\frac{p}{q}$ is a positive rational number.
One may consider a
tubular neighborhood nd$(K)$ of $K$ in $Y$ which may be identified
with $S^1\times D^2$ in such a way that the curve $S^1\times\{1\}$
on the boundary of this solid torus has zero linking number with
$K$ in $Y$ (thus has trivial image in the first homology of
$Y\setminus \text{nd}(K)$). Denote this curve by $\lambda$, and
denote the meridian of $K$- which corresponds to $\{1\}\times
\partial(D^2)$- by $\mu$. Clearly the closed curve $\mu$ bounds a
disk in $Y$. Replacing nd$(K)$ with another solid torus such that
the curve $p\mu+q\lambda$ on the boundary of the solid torus
bounds a disk in the new solid torus produces the
$\frac{p}{q}$-surgery on the knot $K$. We denote the resulting
three-manifold by $Y_{\frac{p}{q}}(K)$. The central circle of this
solid torus is a simple closed curve which would
give a rationally null-homologous knot in $Y_{\frac{p}{q}}(K)$
 denoted by $(Y_{\frac{p}{q}}(K),K_{\frac{p}{q}})$ or just by
$K_{\frac{p}{q}}$ in this paper.\\

The Heegaard Floer homology of $Y_{\frac{p}{q}}(K)$ is computed  in
\cite{OS-Qsurgery} in terms of the Heegaard Floer complex associated
with the knot $(Y,K)$. It is the goal of this paper to understand
the filtration of the complex given in \cite{OS-Qsurgery} induced
by the knot $K_{\frac{p}{q}}$. This is a special case of the question
raised in \cite{filtration}. It will be used in \cite{Ef-gluing}
to study the general case where two knot-complements are glued
along their torus boundary. \\

In order to state the results obtained in this paper, let us introduce
some notation. For simplicity assume that $Y$ is a homology sphere.
Let $H=(\Sig,\alphas,\betas,p)$ be a pointed Heegaard
diagram for the knot $(Y,K)$, and let the corresponding Heegaard
complex be generated by $[\x,i,j]\in (\Ta\cap\Tb)\times \Z\times \Z$.
Assume that the complex $\CFKT^\infty(Y,K)$ is equipped with
the differential $\partial^\infty$. For any intersection point
$\x\in \Ta\cap\Tb$, let $i(\x)\in \Z\cong \RelSpinC(Y,K)$ denote
the associated relative $\SpinC$ structure.\\

Construct a complex $\D$ using generators $[\x,i,j,k,l]\in
(\Ta\cap\Tb)\times \Z^4$ equipped with a differential $\partial_\D$
which is defined by
$$\partial_{\D}[\x,i,j,k,l]=\sum_{p}
n_p[\y_p,i-i_p,j-i_p,k-k_p,l-k_p],$$
whenever
$$\partial^\infty [\x,i,k]=\sum_{p}n_p[\y_p,i-i_p,k-k_p].$$
For $a=[\x,i,j,k,l]\in \D$ define $\Delta(a)=i-j+k-l$ and let
$\D_\delta$ be the subcomplex of $\D$ generated by those $a\in \D$
as above such that $\Delta(a)=\delta$. Define $\D^{up}$ and $\D^{down}$
to be copies of the complex $\D_1$.\\

 Define a
$\Z\oplus\Z$ grading on the generators $\alpha=[\x,i,j,k,l]$
of $\D^{up}\cup \D^{down}$ by
\begin{displaymath}
\begin{split}
&\Gi(\alpha)=
(\text{\emph{max}}(i,l),\text{\emph{max}}(j,k)),\ \ \text{if }
  \alpha\in \D^{up},\\
&\Gi(\alpha)=(i,j)\ \ \ \ \text{if } \alpha\in\D^{down}.
\end{split}
\end{displaymath}

Suppose that $u$ and $w$ are points on the two sides of $p$ which
give the filtration on $\CFKT^\infty(Y,K)$. The $\Z\oplus\Z$
filtration on $\D^{down}$ uses the first marked point $u$. If we
denote by $\D_0^{down}$ the same chain complex with the
$\Z\oplus\Z$ filtration coming from the last two integer
components of the generators (corresponding to the marked point
$w$) then there is a homotopy
 equivalence of filtered chain complexes
$$\tau: \D_0^{down}\lra \D^{down},$$
coming from the invariance of the chain homotopy type from the choice
of the marked point.\\

Define the chain maps $g_{\frac{p}{q}}:\D^{up}\otimes
\frac{\Z[\zeta]}{\zeta^q=1} \lra \D^{down}\otimes
\frac{\Z[\zeta]}{\zeta^q=1}$ via the formula
\begin{displaymath}
\begin{split}
&g_{\frac{p}{q}}([\x,i,j,k,l]\otimes\zeta^t)=
\tau[\x,l,k,2k-j-\lfloor\frac{t+p}{q}\rfloor,2l-i-
\lfloor\frac{t+p}{q}\rfloor]\otimes \zeta^{t+p}.\\
\end{split}
\end{displaymath}
Let $\fbar_{\frac{p}{q}}=Id+g_{\frac{p}{q}}$
be the sum of this map with the identity, and let
$\M(\fbar_{\frac{p}{q}})$ denote the mapping cone of
$\fbar_{\frac{p}{q}}$. The filtration $\Gi$ may naturally be
extended to $\M(\fbar_{\frac{p}{q}})$.

The complex $\M(\fbar_{\frac{p}{q}})$ is
 decomposed as a direct sum according to
the relative $\SpinC$ structures. The $\SpinC$ classes are assigned
to the generators of $\M(\fbar_{\frac{p}{q}})$ via
\begin{displaymath}
\begin{split}
&\relspinc([\x,i,j,k,l]\otimes \zeta^t)=q(i(\x)+j-k)+p(i-j)+t\\
&\text{for }[\x,i,j,k,l]\in \D^{up} \text{ or } \D^{down}, \ \ \
0\leq t<q.
\end{split}
\end{displaymath}

Let $P$ denote a domain in the lattice $\Z\oplus\Z$ such that
for any $p=(i,j)\in P$ there exists some $N_p>0$ such that
for any $i',j'>N_p$ the lattice point $(i-i',j-j')$ is not in
$P$. Assume furthermore that if $(i,j)$ and $(i',j')$ are
points in $P$ then for any pair of integers $i''\in [i,i']$ and
$j''\in [j,j']$ the point $(i'',j'')$ is also included
in $P$. Such a domain $P$ will be
called a \emph{positive test domain}. If $C$ is a $\Z\oplus
\Z$-filtered chain complex with filtration $\mathcal{F}$, denote by $C^P$
the module generated by $\mathcal{F}^{-1}(P)$, and equipped with
the  induced structure of a chain complex
coming from $C$. We will denote by
$H_*^P(C)$ the homology
of the complex $C^P$. Two $\Z\oplus\Z$-filtered chain
complexes $C$ and $C'$ are called \emph{quasi-isomorphic} via a chain
map $f:C\lra C'$ if for any positive test domain
$P\subset \Z\oplus \Z$ the induced map in homology
$$f_*^P:H_*^P(C)\lra H_*^P(C')$$
is an isomorphism.
\\

Among other things, the main theorem
proved in this paper may be stated as follows:

\begin{thm}
Let $Y$ be a homology sphere and let $(Y,K)$ denote a knot in
$K$. Suppose that $\frac{p}{q}> 0$ is a rational number and let
$(Y_{\frac{p}{q}}(K),K_{\frac{p}{q}})$, $\D^{up}$, $\D^{down}$ and
$\M(\fbar_{\frac{p}{q}})$ be as before.
Then  for any relative $\SpinC$ structure
$$\relspinct\in\RelSpinC(Y_{\frac{p}{q}}(K),K_{\frac{p}{q}})
=\RelSpinC(Y,K)\cong \Z$$ the knot Floer complex $\CFKTT^\infty(
Y_{\frac{p}{q}}(K),K_{\frac{p}{q}},\relspinct)$ is quasi-isomorphic
to the complex $\M(\fbar_{\frac{p}{q}})[\relspinct]$.
\end{thm}

In particular, when $P=\{(0,0)\}$ the above theorem is used in the
final section to prove the following non-vanishing result about
$\widehat{\HFKT}$ of rational surgeries on a knot:\\

\begin{thm}
Suppose that $K$ is a knot in $S^3$ of genus $g(K)$, and let
$r=\frac{p}{q}\in \Q$ be a positive rational number. Under the
natural identification
$\RelSpinC(S^3,K_r)=\Z$ we will have
$$\widehat{\HFKTT}(K_r,-qg(K))\cong \widehat{\HFKTT}(K_r,qg(K)+p-1)
\cong \widehat{\HFKTT}(K,g(K))\neq 0,$$
and for any $\relspinct\in \Z$ such that $\relspinct<-qg(K)$
or $\relspinct\geq qg(K)+p$ we will have
$\widehat{\HFKTT}(K_r,\relspinct)=0$.\\
\end{thm}

We may use this theorem to give an easy proof of  a result of
Gordon and Luecke (\cite{GL}) than $\frac{1}{q}$-surgery on
non-trivial knots in $S^3$ can note produce $S^3$. This is a
special case of \emph{Property P} proved by Kronheimer and Mrowka
\cite{KM-PP}. In this form, it is also proved in
\cite{OS-Qsurgery} using Floer homology. However our proof is
different from theirs.\\

\begin{cor}
If $K$ is a knot in $S^3$ and $r=\frac{p}{q}\in \Q$ is a positive
rational number such that the three-manifold obtained by $r$-surgery
on $K$ is $S^3$ then $K$ is the unknot and $p=1$.
\end{cor}
\begin{proof}
By considering the first Betti number it is clear that $p=1$. If the
three-manifold $S^3_r(K)$ is $S^3$, it is implied that
$K_{\frac{1}{q}}$ is a knot $L$ in $S^3$, and for any positive integer
$q'$, $K_{\frac{1}{q+q'}}=L_{\frac{1}{q'}}$. The previous
non-vanishing theorem implies that
$(q+q')g(K)=q'g(L)$ for any $q'$, which implies that $g(K)=0$.
\end{proof}

We first consider the case of large integral surgeries
in section 2.
Together with a surgery exact sequence, this would suggest a
computation of Heegaard Floer complex for arbitrary integral
surgeries on null-homologous knots, which is done in section 3.
In section 4 we will generalize this
formula to the case of rationally null-homologous knots with minor
modifications.
In section 5 we combine these results to prove the above theorem
for rational surgeries. In section 6 we study $\widehat{\HFKT}(K_r)$
for positive rational numbers $r$ and will prove the above
non-vanishing result.
\\

The formulas obtained here are essential for the study of the structure
of Heegaard Floer complex that is associated with the closed
three-manifold obtained by gluing
two three-manifolds along their torus boundaries, which will appear
in the sequel \cite{Ef-gluing}.\\

\textbf{Acknowledgement.} The author would like to thank Zolt\'an
Szab\'o for helpful discussions.

\section{Heegaard Floer homology of large integral surgeries}
The first step toward the desired computation is an understanding of
the case of $n$ surgeries, when $n$ is a large integer.
Suppose that $(Y,K)$ is as above
and consider a Heegaard diagram for the pair. Suppose that
the curve $\beta_g$ in the Heegaard diagram
$$H=(\Sig,\alphas=\{\alpha_1,...,\alpha_g\},\betas=\{\beta_1,...,
\beta_g\},p)$$
corresponds to the meridian of $K$ and that the marked point $p$ is
placed on $\beta_g$. One may assume that the curve $\beta_g$ cuts
$\alpha_g$ once and that this is the only element of $\alphas$ that has
an intersection point with $\beta_g$. Suppose that $\lambda$
represents a longitude for the knot $K$ (i.e. it cuts $\beta_g$ once
and stays disjoint from other elements of $\betas$) such that the
Heegaard diagram
$$(\Sig,\alphas,\{\beta_1,...,\beta_{g-1},\lambda\})$$
represents the three-manifold $Y_0(K)$. Winding $\lambda$ around
$\beta_g$- if it is done $n$ times- would produce a Heegaard diagram
for the three-manifold $Y_n(K)$. More precisely, if the resulting
curve is denoted by $\lambda_n$, the Heegaard diagram
$$H_n=(\Sig,\alphas,\betas_n=\{\beta_1,...,\beta_{g-1},\lambda_n\},p_n)$$
would give a diagram associated with the knot $(Y_n(K),K_n)$, where
$p_n$ is a marked point placed on $\lambda_n$.\\

Denote by $u_n$ and $w_n$ a pair of marked points on the two sides of
$\lambda_n$ which are both very close to $p_n$. In the description
of the Heegaard complex associated with $Y_n(K)$ given in
\cite{OS-knot} we are interested in computing the
$\Z\oplus \Z$-filtration induced by the pair of points $(u_n,w_n)$.\\

Fix a $\SpinC$-structure $\spinc\in \SpinC(Y_n(K))$ and choose the
marked point $p_n$ so that all the generators of the complex
$\CFT^\infty(Y_n(K),\spinc)$ are supported in the winding region
if the $\SpinC$ classes are assigned using
either of the marked points $u_n$ or $w_n$ on the two sides of
$\lambda_n$.\\

The curve $\lambda_n$ intersects the $\alpha$-curve $\alpha_g$ in
$n$-points which appear in the winding region (there may be other
intersections outside the winding region). Denote these points
of intersection by
$$...,x_{-2},x_{-1},x_{0},x_{1},x_{2},...,$$
where $x_1$ is the intersection point with the property that three of
its four neighboring quadrants belong to the regions that
 contain either $u_n$ or $w_n$. Any
generator which is supported in the winding region is of the form
$$\{x_i\}\cup \y_0=\{y_1,...,y_{g-1},x_i\},$$
and it is in correspondence with the generator
$$\y=\{x\}\cup\y_0=\{y_1,...,y_{g-1},x\}$$
for the complex associated with the knot $(Y,K)$, where $x$ denotes the
unique intersection point of $\alpha_g$ and $\beta_g$. Denote the former
generator by $(\y)_i$, keeping track of the intersection point $x_i$
among those in the winding region.\\

Note that $\RelSpinC(Y,K)=\SpinC(Y_0(K))=\Z\oplus \SpinC(Y)$ and for
$\relspinc \in \RelSpinC(Y,K)$ let $i(\relspinc)$ denote the first
component in this decomposition. Remember that there is a map
$$\relspinc:\Ta\cap \Tb \lra \RelSpinC(Y,K),$$
which is defined in
\cite{OS-knot} or more generally in \cite{OS-Qsurgery}.
Thus to any generator $\y$ as above we may associate
an integer $i(\relspinc(\y))\in \Z$.\\

Similarly, $\SpinC(Y_n(K))=\frac{\Z}{n\Z}\oplus \SpinC(Y)$ and for
$\spinc_n\in \SpinC(Y_n(K))$ we may denote the projection over the
first component of this decomposition by $i_n(\spinc_n))\in
\frac{\Z}{n\Z}$.  The marked point gives a map
$$\spinc_n:\Ta\cap\T_{\beta_n}\lra \SpinC(Y_n(K)).$$
As a result, for any generator $\y$ as above and any integer $i\in
\Z$ which is not \emph{very large} we obtain a number
$$i_n(\spinc_n((\y)_i))=[i_n(\spinc_n((\y)_0))-i] \in
\frac{\Z}{n\Z}.$$

Note that there is a relation between these numbers given by the
following lemma which is in fact proved in \cite{OS-knot}:

\begin{lem}
With the above notation, for any generator $\y \in \Ta\cap \Tb$ we have
$$i(\relspinc(\y))-i=i_n(\spinc_n((\y)_i))\ \  (\text{\emph{mod }}n).$$
\end{lem}

Let $s:\RelSpinC(Y,K)\lra \SpinC(Y)$ and $s_n:\RelSpinC(Y_n(K),K_n)\lra
\SpinC(Y_n(K))$ denote the natural maps obtained by extending the
relative $\SpinC$ structures of the knot complements over the attached
solid torus.\\

For a fixed $\SpinC$-structure $\spinct_n\in \SpinC(Y_n(K))$ define
\begin{displaymath}
\begin{split}
\CFKT^\infty&(Y,K,\spinct_n):=\\
&\Big\langle\Big\{[\y,i,j]\in (\Ta\cap \Tb)\times
\Z \times \Z\ \Big|\ \ \begin{array}{c}
s_n(\relspinc(\y)-(i-j)\text{PD}[\mu])=\spinct_n\\
-\frac{n}{2}\leq i(\relspinc(\y))-i+j<\frac{n}{2}\\
\end{array} \Big\}\Big\rangle,
\end{split}
\end{displaymath}
where $\mu$ denotes the image of the meridian of the knot $K$ in the
three-manifold $Y_0(K)$ obtained by a zero-surgery on $K$.\\

It is shown in \cite{OS-knot} that for large values of $n$
the complex
$\CFKT^\infty(Y,K,\spinct_n)$ may in fact be thought
of as giving the complex
$\CFT^\infty (Y_n(K),\spinct_n)$ (as a $\Z$-filtered chain
complex)
under a correspondence which may be naively described as
$$[\y,i,j]\mapsto [(\y)_{i-j},\ \text{max}\{i,j\}].$$

This map is constructed by counting holomorphic triangles
corresponding to the triple Heegaard diagram
$$R_n=(\Sig, \alphas,\betas,\betas_n;u',w').$$
Here $\betas_n=\{\beta_1',...,\beta_{g-1}',\lambda_n\}$
denotes a set of $g$ simple closed curves such that
the first $g-1$ of them are in fact
close isotopic copies of the curves in
$\betas$ with the property that each
$\beta_i'$ cuts the corresponding curve $\beta_i$ in a pair of cancelling
intersection points, and $u',w'$ are two marked points on the two
sides of the curve $\beta_g$ very close to the point $p\in \beta_g$.
The curve $\beta_g$ is chosen so that it is located almost in the
middle of the winding region. We may abuse the notation and
write $$\betas_n=\{\beta_1,...,\beta_{g-1},\lambda_n\}.$$
The inverse of the above map is in fact of the form
$$[\y,i,j]\mapsto [(\y)_{i-j},\ \text{max}\{i,j\}]+
\ \text{lower order terms},$$ which would then imply that the
structure of a filtered chain complex induced on the
right-hand-side using the initial map $[\y,i,j]\mapsto
[(\y)_{i-j},\ \text{min}\{i,j\}]$ is the same as its own structure
as a filtered chain complex, see \cite{OS-knot} for a more
detailed description.\\

\begin{figure}
\mbox{\vbox{\epsfbox{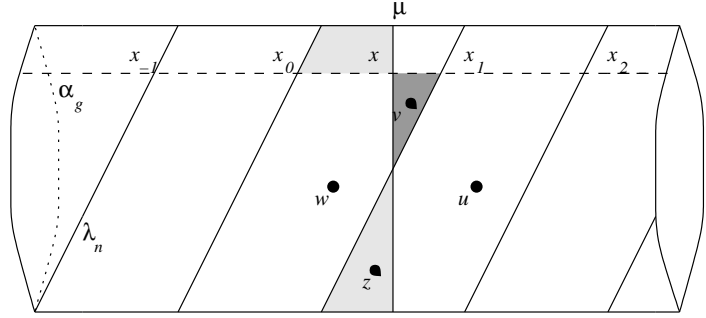}}} \caption{\label{fig:HD-Rn} {
The Heegaard diagram $S_n$. The shaded triangles are $\Delta_0$
and $\Delta_1$.
}}
\end{figure}

Choose the curve $\beta_g$
so that it cuts $\lambda_n$ exactly once in $\{y\}=\beta_g\cap
\lambda_n$ and the intersection point $\{x\}=\alpha_g\cap
\beta_g$ is located between $x_0$ and $x_1$. Furthermore, there is a small
triangle $\Delta_0$ with vertices $x,x_0,y$ on the surface $\Sig$,
and another one- denoted by $\Delta_1$- with vertices $x,x_1,y$.
There are $4$ quadrants around $y$, two of them being parts of the
triangles $\Delta_0$ and $\Delta_1$.
Denote the remaining two regions by $D_1$ and $D_2$,
so that $D_1$ is on the
right-hand-side of both $\beta_g$ and $\lambda_n$ and so that $D_2$
is on the left-hand-side of both of them. Let $u,v,w$ and $z$ be four
marked points in $D_1,\Delta_1,D_2$ and $\Delta_0$ respectively (
see figure~\ref{fig:HD-Rn}).

Denote the new
Heegaard diagram by
$$S_n=(\Sig,\alphas,\betas_n,\betas;u,v,w,z).$$
There is a holomorphic triangle map associated with $S_n$ which
we will study below. Note that the Heegaard diagram
$$\overline{H}=(\Sig,\alphas,\betas;u,v,w,z)$$
corresponds to a chain complex $\D$. The generators of $\D$ are of the
form $$[\x,i,j,k,l]\in (\Ta\cap \Tb)\times \Z\times \Z \times \Z\times
\Z=(\Ta\cap\Tb)\times \Z^4.$$
If $\partial^\infty$ denotes the boundary map associated with the
filtered chain complex $\CFKT^\infty(Y,K)$
assigned to the Heegaard diagram $H$ for the
knot and
$$\partial^\infty [\x,i,k]=\sum_{p}n_p[\y_p,i-i_p,k-k_p],$$
then the differential $\partial_{\D}$ for $\D$ is defined by
$$\partial_{\D}[\x,i,j,k,l]=\sum_{p}n_p[\y_p,i-i_p,j-i_p,k-k_p,l-k_p].$$
We may note that if two generators $[\x,i,j,k,l]$ and
$[\y,i',j',k',l']$ are connected by a topological disk $\phi
\in \pi_2(\x,\y)$ then we will have the identities
\begin{equation}
\label{eq:Didentity}
\begin{cases}
\spinc(\x)=\spinc(\y)\\
i(\relspinc(\x))-(i-k)=i(\relspinc(\y))-(i'-k')\\
i-j=i'-j',\ \ k-l=k'-l',
\end{cases}
\end{equation}
where we define $\spinc(\x)=s(\relspinc(\x))$.\\

In $\D$ let $\D_1$ denote the subcomplex generated by all
$a=[\x,i,j,k,l]$ such that $\Delta(a)=i-j+k-l$ is equal to
$1$. The triangle map $\Phi$ maps  a third complex
$\E_0\subset \E$ to the complex $\D_1$. Here $\E$
is the filtered chain complex associated with the Heegaard
diagram
$$H''=(\Sig,\alphas,\betas_n;u,v,w,z),$$
which is a Heegaard diagram associated to
$\CFKT^\infty(Y_n(K),K_n)$ in a way similar to the above
correspondence between $\D$ and $\CFKT^\infty(Y,K)$. Here the two
marked points $w$ and $v$ are in the same domain and similarly $u$
and $z$ are in the same domain (the last condition above is thus
replaced by $j-k=j'-k'$ and $i-l=i'-l'$). Denote by $\E_0$ the
subcomplex of $\E$ consisting of the tuples
$$\E_0=\Big\langle\Big\{
[\x,i,j,k,l]\in \E\ \Big|\ \begin{array}{c} j-k=0\\ i-l=0
\end{array}
\Big\}\Big\rangle.$$
The complex $\E_0$ may be identified with the chain
complex $\CFKT^\infty(Y_n(K),K_n)$. We will sometimes abuse the
notation and denote the element $[\x,i,j,j,i]\in \E_0$ by $[\x,i,j]$. \\

The triangle map $\Phi:\E\lra \D$ reduces to a map $\Phi_0:\E_0\lra
\D_1$. This may be checked by examining the local multiplicities around
the intersection of $\mu=\beta_g$ and $\lambda_n$.\\

In terms of the Energy filtration and the $\Z\oplus\Z$-filtration,
the map $\Phi_0$ is defined so that:
\begin{displaymath}
\begin{split}
\Phi_0^{-1}[\x,i,j,k,l]=[(\x)_{k-j},\text{max}\{i,l\},\ &\text{max}
\{j,k\}]+\\ &\text{terms of lower order,}
\end{split}
\end{displaymath}
when $|k-j|$ is small. Having this in mind, equip $\D_1$ with a
$\Z\oplus\Z$ filtration
$$\mathcal{F}[\x,i,j,k,l]=(\text{max}\{i,l\},\ \text{max}\{j,k\}).$$

The set of relative $\SpinC$ structures for the knot $(Y_n(K),K_n)$ is
easy to understand, according to \cite{OS-Qsurgery}:
\begin{equation}
\begin{split}
\RelSpinC(Y_n(K),K_n)&=\RelSpinC(Y_n(K)\setminus \text{nd}(K_n),
\partial(Y_n(K)\setminus \text{nd}(K_n)))\\
&=\RelSpinC(Y\setminus
\text{nd}(K),\partial (Y\setminus \text{nd}(K)))\\
&=\RelSpinC(Y,K)=\SpinC(Y)\oplus \Z
\end{split}
\end{equation}
There is a map which projects
relative $\SpinC$ structures over $\SpinC$ structures of $Y_n(K)$:
\begin{displaymath}
\begin{split}
G_n=G_{Y_n(K),K_n}(=s_n):\SpinC(Y)&\oplus \Z=\RelSpinC(Y_n(K),K_n)\\
&\lra \SpinC(Y_n(K))=\SpinC(Y)\oplus \frac{\Z}{n\Z}
\end{split}
\end{displaymath}
The map $G_n$ is simply the reduction modulo the integer $n$, i.e.
$G_n(\spinc,i)=(\spinc,(i)_{\text{mod}\ n})$. There is a relative
$\SpinC$ structure associated with any generator of the
complex $\CFKT^\infty(Y_n(K),K_n)$.  The
intersection point associated with $[\x,i,j,k,l]$ is $(\x)_{k-j}$ and
the relative $\SpinC$ structure associated with this generator
is $$\relspinc(\x)+(j-k+n\epsilon(k-j))
\text{PD}[\mu]\in \RelSpinC(Y,K),$$ where $\epsilon(i)$ is equal to
$1$ if $i$ is positive and is zero otherwise. This may be checked
easily from lemma 2.1 in \cite{OS-knot}, at least in the relative
version.\\

 We may assign relative $\SpinC$ structures to the generators of
 $\D_1$ so that the map $\Phi$ preserves the relative $\SpinC$ class.
With the above computation and the filtration $\mathcal{F}$ in mind,
this assignment should be defined via
\begin{displaymath}
\begin{split}
\relspinc:(\Ta\cap\Tb)\times \Z^4 \lra &\RelSpinC(Y_n(K),K_n)=
\SpinC(Y)\oplus\Z\\
\relspinc[\x,i,j,k,l]=\relspinc(\x)&+(j-k+n\epsilon(k-j))\text{PD}[\mu]\\
&+n(\text{max}\{i,l\}-\text{max}\{j,k\})\text{PD}[\mu]
\end{split}
\end{displaymath}
Note that the last term is $n(i-j-\epsilon(k-j))$ since $i-j+k-l=1$.
It is implied that
\begin{displaymath}
\relspinc[\x,i,j,k,l]=\relspinc(\x)+(j-k+n(i-j))\text{PD}[\mu].
\end{displaymath}

The $\SpinC$ structure in $\RelSpinC(Y_n(K))$
associated with a generator $a=[\x,i,j,k,l]$ as
above is
$$\spinc_n(a)=(s(\x),(i(\relspinc(\x))+j-k)_{\text{mod }n})\in
\SpinC(Y)\oplus \frac{\Z}{n\Z}=\SpinC(Y_n(K)).$$

For $s\in \Z$ define $\D_1^s$ to be the subcomplex of $\D_1$ generated
by $[\x,i,j,k,l]$ such that $i(\relspinc(\x))+j-k=s$. Fix a relative
$\SpinC$ structure $\relspinc\in \RelSpinC(Y_n(K),K_n)=\RelSpinC(Y,K)$
and let $-\frac{n}{2}\leq s<\frac{n}{2}$ be an integer with the
property that $i(\relspinc)=s\ (\text{mod }n)$. Compose the map
$\Phi:\CFKT^\infty(Y_n(K),K_n,\relspinc)\rightarrow \D_1$ with the
projection over $\D_1^s(\relspinc)$.
Here $\D_1^s(\relspinc)$ is generated by those generators of
$\D_1^s$ which are in the relative $\SpinC$ class $\relspinc$, i.e.
with $[\x,i,j,k,l]$ such that
$\relspinc(\x)+(j-k+n(i-j))\text{PD}[\mu]=
\relspinc$.
It is not hard to show that this induces a
chain homotopy equivalence
$$\Phi^{\relspinc}:\CFKT^\infty(Y_n(K),K_n,\relspinc)\lra \D_1^s(\relspinc).$$
We have proved the following theorem.


\begin{thm}
\label{thm:Large-Surgery}
Suppose that $(Y,K)$ is a null-homologous knot, $n$ is a
sufficiently large integer, and $(Y_n(K),K_n)$,
the Heegaard diagram $H$ and the
complexes $\D_1$ and $\D_1^s(\relspinc)$ are as above.
Then the knot Floer complex
associated with the rationally null-homologous knot $(Y_n(K),K_n)$
in the relative $\SpinC$ class $\relspinc\in \RelSpinC(Y_n(K),K_n)$
has the same filtered chain homotopy type as the complex $\D_1^s(\relspinc)$
equipped
with the $\Z\oplus\Z$ filtration given by
$$\mathcal{F}[\x,i,j,k,l]=(\text{max}\{i,l\},\
\text{max}\{j,k\}).$$
Here $-\frac{n}{2}\leq s<\frac{n}{2}$ is chosen so that
$i(\relspinc)=s\ (\text{\emph{mod }}n)$ and the decomposition of
$\D_1$ into a direct sum according to the relative $\SpinC$ structures
is given by
$$\D_1(\relspinct)=\Big\langle\Big\{[\x,i,j,k,l]\in \D_1\ \Big| \
\relspinc(\x)+(j-k+n(i-j))\text{\emph{PD}}[\mu]=\relspinct
\Big\}\Big\rangle.$$

\end{thm}

This computation finishes our study of the surgery formulas when the
integer $n$ is large. We will use this computation in the upcoming
section to understand the chain complex associated with the knot
$(Y_n(K),K_n)$, where $n$ is an arbitrary non-zero integer.\\

\begin{remark}
Note that the computation of Ozsv\'ath and Szab\'o in \cite{OS-knot}
of the complex
$\text{\emph{CF}}^+(Y_n(K))$ is in fact a corollary of the above theorem, once
we write down the relation between the filtered chain complex
$\CFKTT^\infty(Y_n(K),K_n)$ and the Heegaard Floer homology associated
with the ambient three-manifold $Y_n(K)$ as introduced in \cite{OS-Qsurgery}.
\end{remark}

\begin{remark}
If $K$ is an alternating knot in $S^3$, the knot Floer complex
associated with $K$ may be computed as in \cite{Ras, OS-alternating}.
As a result the above argument gives a computation of the complexes
$\CFKTT^\infty(S^3_n(K),K_n)$ where $n$ is sufficiently large, which is
not hard to do by hand once the symmetrized Alexander polynomial
and the signature are
given.
\end{remark}

\section{Modification of Ozsv\'ath-Szab\'o argument}
We begin by reminding the reader of the argument given by Ozsv\'ath
and Szab\'o in \cite{OS-Zsurgery} to generalize the computation of
$\text{CF}^+(Y_n(K))$ when $n$ is a large integer to the case of arbitrary
integer $n$. \\

In \cite{OS-Zsurgery}, the first step is the construction of an
exact sequence coming from the quadruple Heegaard diagram
$$R_{m,n}=(\Sig, \alphas, \betas, \betas_n,\betas_{m+n})$$
where $\betas, \betas_n$ and $\betas_{n+m}$ are as before, and
each of them is equipped with a special curve denoted by
$\mu=\beta_g,\lambda_n$ and $\lambda_{m+n}$ respectively.
Furthermore, we assume that the curves in $\betas_n$ are small
isotopic copies of the curves in $\betas$ (except for $\lambda_n$
which is not an isotopic copy of $\mu$) such that there is a pair
of cancelling intersection points between any curve in $\betas_n$
and the corresponding element of $\betas$. The same assumption is
made for $\betas_{m+n}$ and the same is assumed for the relation
between the curves in $\betas_n$ and those in $\betas_{m+n}$. A
marked point $p$ is fixed on $\mu$ and a marked point $u$ on the
diagram is also chosen. Using the marked points $u$ and $p$ three chain
maps are constructed:
\begin{displaymath}
\begin{split}
&f_1^+:\CFT^+ (Y_n(K))\lra \CFT^+(Y_{m+n}(K))\\
&f_2^+:\CFT^+ (Y_{m+n}(K))\lra \bigoplus^m \CFT^+(Y)\\
&f_3^+:\bigoplus^m \CFT^+(Y)\lra \CFT^+(Y_n(K)),
\end{split}
\end{displaymath}
where the map $f_2^+$ counts the number of times the boundary
of a holomorphic triangle touches the point $p$, in order to give
a map to the chain complex with twisted coefficients
$$\CFT^+(Y,\Z[\frac{\Z}{m\Z}])\cong\CFT^+(Y)\otimes_\Z \frac{\Z[T]}{T^m=1}.$$

As the integer $m=\ell n$ becomes large the complex $\CFT^+(Y_n(K))$
is proved to be chain homotopic to the limit of the
mapping cones of the chain maps $f_2^+$ via
the map induced by $f_1^+$. Since for large choices of $\ell$ the complex
$\CFT^+(Y_{(\ell+1)n}(K))$ is described in terms of $\CFKT^\infty(Y,K)$ in
\cite{OS-knot}, one would be done with the computation once the
compatibility of certain maps and isomorphisms are verified, and
the limiting behavior as $m$ goes to infinity is studied.\\

\begin{figure}
\mbox{\vbox{\epsfbox{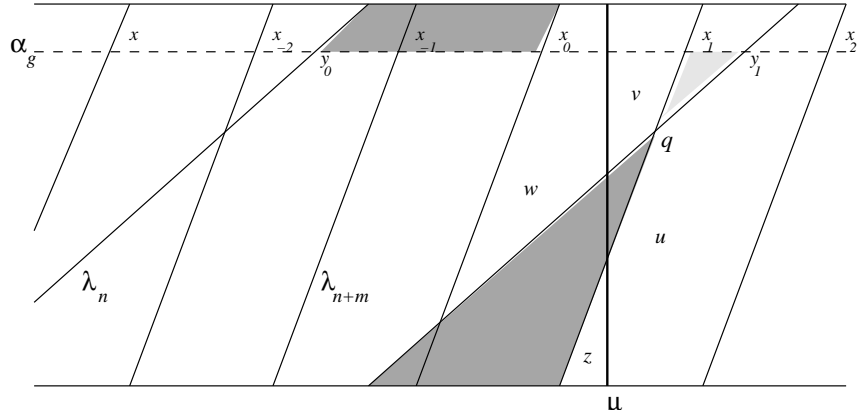}}} \caption{\label{fig:HD-Rmn} {
The Heegaard diagram $R_{m,n}$. The shaded triangles are $\Delta_0$
and $\Delta_1$. The marked points $u,v,w$ and $z$ are placed in
$D_1^R,D_2^R,D_2^L$ and $D_1^L$ respectively.
}}
\end{figure}

In this section we will basically follow the same strategy using more
marked points. Add the marked points $u,v,w$ and $z$ to the Heegaard
diagram $R_{m,n}$ in the regions described
bellow. Choose an intersection point between the curves $\lambda_n$
and $\lambda_{m+n}$ in the middle of the winding region, denoted by
$q$. We will assume that $m=\ell n$ for an integer $\ell$ which is chosen to
be appropriately large. We continue to assume that $\alpha_g$ is the
unique $\alpha$-curve in the winding region.
From the $4$ quadrants around the intersection
point $q$, two of them are parts of  small triangles $\Delta_0$ and
$\Delta_1$ between
$\alphas,\betas_n$ and $\betas_{m+n}$.
We may assume that the intersection points between $\alpha_g$ and
$\lambda_{m+n}$ in the winding region are
$$...,y_{-2},y_{-1},y_0,y_1,y_2,...,$$
and that the intersection points between $\lambda_n$ and $\alpha_g$
are $$...,x_{-2},x_{-1},x_0,x_1,x_2,...$$
as before. We may also assume that the domain
$\Delta_i$ for $i=0,1$ is the triangle with vertices
$q,x_i$ and $y_i$, and $\Delta_1$ is one of the connected
domains in the complement of curves $\Sig\setminus C$ where
 $$C=\alphas \cup\betas_n\cup \betas_{m+n}.$$
Other that $\Delta_0$ and $\Delta_1$ there are two other
domains which have $q$ as a corner. One of them is on the
right-hand-side of both $\lambda_n$ and $\lambda_{m+n}$, denoted by
$D_1$, and the other one is on the left-hand-side of both of them,
denoted by $D_2$. The domains
$D_1$ and $D_2$ are assumed to be connected regions
in the complement of the
curves $\Sig\setminus C$.
We may assume that the meridian $\mu$ passes through the regions
$D_1,D_2$ and $\Delta_0$, cutting each of them into two parts:
$\Delta_0=\Delta_0^R\cup\Delta_0^L$, $D_1=D_1^R\cup D_1^L$ and
$D_2=D_2^R\cup D_2^L$. Here $\Delta_0^R\subset \Delta_0$ is the
part on the right-hand-side of $\mu$ and $\Delta_0^L$ is the part
on the left-hand-side. Similarly for the other partitions. Choose
the marked points so that $u$ is in $D_1^R$, $v$ is in $D_2^R$, $w$
is in $D_2^L$ and $z$ is in $D_1^L$ (see figure~\ref{fig:HD-Rmn}).
We obtain the Heegaard diagram
$$R_{m,n}=(\Sig,\alphas,\betas,\betas_n,\betas_{m+n};u,v,w,z),$$
which will be used for constructing certain relevant chain complexes
and chain maps between them.\\

We are interested in the mapping cone of the chain map $f$ defined
from the filtered chain complex $\B_{m+n}$ associated with the pair
$(\alphas,\betas_{m+n})$ to the chain complex $\B$ associated with
$(\alphas,\betas)$. The map $g$ going from the complex $\B_n$
associated with $(\alphas,\betas_n)$ to $\B_{m+n}$ defines a map
$\overline{g}$ from $\B_n$ to the mapping cone $\M(f)$ of $f$, and the
chain map $h$ from $\B$ to $\B_n$ defines a map $\overline{h}$ from
$\M(f)$ to $\B_n$. More precisely, the Heegaard diagram
$$(\Sig, \alphas,\betas_n;(u,z),(v,w))$$ will produce the
$\Z\oplus\Z$-filtered chain complex $\B_n$ which has the chain
homotopy type of $\CFKT^\infty(Y_n(K),K_n)$. Here by putting the
pairs of points in parenthesis we mean that in the relevant
Heegaard diagram the two points are in the same domain in the
complement of the curves appearing in the Heegaard diagram. The
triangle map which defines $g$ is associated with the Heegaard
triple
$$(\Sig,\alphas,\betas_n,\betas_{m+n};(u,z),(v,w)),$$
and will provide a map from $\B_n$ to the filtered chain complex
$\B_{m+n}$ (defined using the Heegaard diagram
$(\Sig,\alphas,\betas_{m+n};(u,z),(v,w))$).
A canonical generator in the complex associated with $(\Sig,\betas_n,
\betas_{m+n};u,w)$
is chosen as specified in \cite{OS-Qsurgery}, and is used to define
$g$ (this is a Heegaard diagram for the canonical knot $O_{\frac{m}{1}}$ of
the Lens space $L(m,1)$).\\

The Heegaard diagram $(\Sig,\alphas,\betas;u,v)$ defines a $\Z\oplus
\Z$-filtered chain complex with twisted coefficients, with the help of
marked points $w$ and $z$. More precisely,
consider the chain complex $\B[m]$
generated by the generators $[\x,i,j]$ with
$\x\in \Ta\cap\Tb$, and $i,j\in \Z$ over the coefficient ring
$\Z[\frac{\Z}{m\Z}]=\frac{\Z[T]}{T^m=1}$. The boundary map is defined
via
$$\partial^m[\x,i,j]=\sum_{\substack{\y \in \Ta\cap\Tb\\ \phi\in
    \pi_2(\x,\y),\ \mu(\phi)=1\\
   }}
    \#(\widehat{\Mod}(\phi))[\y,i-n_u(\phi),j-n_v(\phi)]T^
{m_p(\phi)},$$
where as usual $n_u(\phi)$ and $n_v(\phi)$ are the intersection
    numbers of the disk $\phi$ with the hyper-surfaces
$\{u\}\times \Sym^{g-1}\Sig$ and $\{v\}\times \Sym^{g-1}\Sig $
    respectively.
Here $m_p(\phi)=n_w(\phi)-n_v(\phi)$ is the intersection number
of $\partial \phi$ with the submanifold $$\beta_1\times \beta_2
\times ...\times \beta_{g-1}\times \{p\}\subset \Tb$$
of $\Tb$. The point $p$ is chosen to be a marked point on
$\mu=\beta_g$ on the arc between $D_2^R$ and $D_2^L$.
Note that $n_u(\phi)=n_v(\phi)$ and that
    $n_w(\phi)-n_v(\phi)=n_z(\phi)-n_u(\phi)$.
The complex $\B[m]$ may be constructed from the knot Floer
complex $\CFKT^\infty(Y,K)$ as the above equations suggest.
\\

Suppose that $\B$ denotes the complex generated over $\Z$ by $[\x,i,j]$
such that $\x\in \Ta\cap\Tb$ and $i,j\in \Z$, equipped with
the differential
$$\partial[\x,i,j]=\sum_{\substack{\y \in \Ta\cap\Tb\\
\phi\in \pi_2(\x,\y)\\ \mu(\phi)=1}}
\#(\widehat{\Mod}(\phi))[\y,i-n_u(\phi),j-n_v(\phi)].
$$
This complex may in fact be constructed out of $\CFT^\infty(Y,K)$.
Then there is a chain map
$$\theta:\B[m]\lra \B\otimes_\Z \Z[\frac{\Z}{m\Z}]
$$
defined by $\theta([\x,i,j]T^k)=[\x,i,j]\otimes T^{k+m_p(\phi_\x)}$,
where $\phi_\x\in \pi_2(\x,\x_0)$ is an arbitrary disk for a
fixed intersection point $\x_0$ of $\Ta$ and $\Tb$ in the $\SpinC$
class of $\x$. Note that $m_p(\phi_\x)$ depends only on $\x,\x_0$ and
is independent of the choice of $\phi_\x\in \pi_2(\x,\x_0)$.\\

The Heegaard triple $(\Sig,\alphas,\betas_{m+n},\betas;u,v,w,z)$
produces a pair of chain map $f_v$ and $f_w$ from the complex
$\B_{m+n}$ to $\B[m]$, defined as follows. The Heegaard diagram
$(\Sig,\betas_{m+n},\betas,(u,v),(w,z))$ is a Heegaard diagram for
the trivial knot in $\#^{g-1}(S^1\times S^2)$. Denote the
canonical generator by $[\Theta_0,0,0]$. Pairing $\Theta_0$ with
the generators of $\B_{m+n}$, using the triangle map associated
with the above Heegaard diagram and the marked points $u,v$, and
twisting according to the difference between the intersection
numbers at $v$ and $w$ we obtain a chain map $f_v$ from $\B_{m+n}$
to $\B[m]$ which is naively defined via
$$f_v[\x,i,j]=\sum_{\substack{\y\in \Ta\cap\Tb\\
\psi\in \pi_2[\x,\Theta_0,\y]\\
\mu(\psi)=0}}\#(\Mod(\psi))[\y,i-n_u(\psi),
j-n_v(\psi)]T^{m_w(\psi)-n_v(\psi)}.$$
The second map $f_w$ is defined similarly:
$$f_w[\x,i,j]=\sum_{\substack{\y\in \Ta\cap\Tb\\
\psi\in \pi_2[\x,\Theta_0,\y]\\
\mu(\psi)=0}}\#(\Mod(\psi))[\y,i-n_z(\psi),
j-n_w(\psi)]T^{m_w(\psi)-n_v(\psi)}.$$

Suppose that the marked point $p\in \mu$ is chosen so that it is
located on the arc between $D_2^R$ and $D_2^L$.
The map $m_p:\pi_2(\x,\Theta_0,\y)\lra \frac{\Z}{m\Z}$ lifts
to a function $$\m:\SpinC(W_{m+n}(K))\lra \Z,$$
where $W_{m+n}(K):Y_{m+n}(K)\lra Y$ is the $4$-manifold
cobordism between these two three-manifolds which is given
by our Heegaard diagram. If $\widehat{F}\subset W_{m+n}(K)$
is the surface obtained by capping a Seifert surface
associated with the knot $K$ in $W_{m+n}(K)$,
the map $\m$ satisfies
the relation
$$\m(\spinct-\text{PD}[\widehat{F}])=\m(\spinct)+m+n.$$
Thus the map $\theta \circ f_v$ may be written as
$$\theta \circ f_v(a)=\sum_{\spinct\in \SpinC(W_{m+n}(K))}
f_v^\spinct(a)\otimes T^{\m(\spinct)},$$
where $f_v^\spinct$ corresponds to counting holomorphic triangles
which induce the $\SpinC$ structure $\relspinct$
on the cobordism.
Similarly we have
$$\theta \circ f_w(a)=\sum_{\spinct\in \SpinC(W_{m+n}(K))}
f_w^\spinct(a)\otimes T^{\m(\spinct)}.$$
Fix a $\SpinC$ structure $\spinct\in \SpinC(Y_{m+n}(K))$ and let
$-\frac{m+n}{2}\leq s<\frac{m+n}{2}$ be an integer such that its
class modulo $m+n$ represents $i_{m+n}(\spinct)\in \frac{\Z}{(m+n)\Z}$.
Denote by $\xfrak_{\spinct}$
and $\yfrak_{\spinct}$ the $\SpinC$-structures on $W_{m+n}(K)$ such
that they restrict to $\spinct$ on $Y_{m+n}(K)$ and
$$\langle c_1(\xfrak_{\spinct}),\text{PD}[\widehat{F}]\rangle
+n+m=2s\ \ \text{and } \
\langle c_1(\yfrak_{\spinct}), \text{PD}[\widehat{F}]\rangle
-n-m=2s.$$
Since for $m=n\ell$ we have $[\widehat{F}].[\widehat{F}]=-n(\ell+1)$, we
obtain the relation $\xfrak_{\spinct}+\text{PD}[\widehat{F}]
=\yfrak_{\spinct}$.
This implies that $$\m(\xfrak_{\spinct})-\m(\yfrak_{\spinct})=-n(\ell+1)=-n\
(\text{mod}\ n\ell).$$
 Denote by
$f_{hor}$ the sum
$$f_{hor.}=\sum_{\spinct\in \SpinC(Y_n(K))}
f_v^{\xfrak_{\spinct}}.T^{\m(\xfrak_{\spinct})}.$$
Similarly, define $f_{ver.}$ using $f_w$, or alternatively
using $f_v$ and the $\SpinC$ structures $\{\yfrak_{\spinct}\}_{\spinct}$.
Denote by $f:\B_{m+n}\lra \B\otimes_\Z \Z[\Z/m\Z]$ the sum of these two maps.
We would like to study the map $f$ under the identification
of $\B_{m+n}(\relspinct)$ with the complex $\D_1^s(\relspinct)$ given by
theorem~\ref{thm:Large-Surgery}.
\\

In fact, in the
proof of theorem~\ref{thm:Large-Surgery} we may use the Heegaard
diagram $R_{m,n}$. Then the composition of the map $f_{hor.}$ with
$\Phi_0^{-1}$ may be described (after composing with another
filtered chain homotopy) as the map
\begin{displaymath}
\begin{split}
&h:\D_1^s\lra \B\otimes_\Z\Z[\frac{\Z}{m\Z}]\\
&h[\x,i,j,k,l]=[\x,i,j]\otimes T^{-s},
\end{split}
\end{displaymath}
and the composition of $f_{ver.}$ may be described as the map
given by
\begin{displaymath}
\begin{split}
&v:\D_1^s\lra \B\otimes_\Z\Z[\frac{\Z}{m\Z}]\\
&v[\x,i,j,k,l]=\tau[\x,l,k]\otimes T^{-s-n}.
\end{split}
\end{displaymath}
Here $\tau$ denotes the chain homotopy equivalence from
$\B$ to itself which comes from changing the pair of marked
points $(u,v)$ to $(z,w)$. In fact defining $f_{ver}$ using
$f_w$ and then composing with $\tau$ is the same as defining
it using $f_v$ and the $\SpinC$ structures $\{\yfrak_{\spinct}
\}_\spinct$ over the cobordism $W_\ell$.\\

Let $\fbar_\ell$ denote the sum of the two maps
$$h,v:
\bigoplus_{-\frac{n(\ell+1)}{2}\leq s <\frac{n(\ell+1)}{2}}
\D_1^s \lra \B\otimes\Z[\Z/m\Z],$$
and denote by $\M(\fbar_\ell)$ the mapping
cone of $\fbar_\ell$. The filtration $\mathcal{F}$ may be extended to
a filtration on $\M(\fbar_\ell)$ defined by
\begin{displaymath}
\begin{split}
&\Gi:(\text{Generators of }\M(\fbar_\ell))\lra \Z\oplus\Z\\
&\begin{cases}
  \Gi([\x,i,j,k,l])=(\text{max}(i,l),\text{max}(j,k))
\ \ \ \ \ \ \text{if }[\x,i,j,k,l]\in \D_1^s,\\
\Gi([\x,i,j]\otimes T^k)=(i,j)\ \ \ \ \ \ \ \ \ \ \
\text{if }[\x,i,j]\otimes T^k\in \B\otimes_\Z\Z[\frac{\Z}{m\Z}].
\end{cases}
\end{split}
\end{displaymath}
The complex $\M(\fbar_\ell)$ decomposes according to the class
of $s$ modulo the integer $n$ and we get
$$\M(\fbar_\ell)=\bigoplus_{[s]\in \frac{\Z}{n\Z}}\M^{[s]}(\fbar_\ell).$$

We would like to identify the decomposition of $\M^{[s]}(\fbar_\ell)$
according to the relative $\SpinC$ structures $\relspinct\in
\RelSpinC(Y_n(K),K_n)=\RelSpinC(Y,K)$. To this end, note that
the image of a generator $[(\x)_l,i,j]\in \B_n$ under the
chain map $g$ is of the form ``$[(\x)_l,i,j]+$ lower order terms''
as an element of $\B_{m+n}$. In particular, the relative $\SpinC$
class of $[\x,i,j,k,l]\in \D_1^s\subset \M^{[s]}(\fbar_\ell)$ under the map
$\Phi_0^{-1}$ should be declared equal to the relative $\SpinC$
class of $[(\x)_{k-j},\text{max}(i,l),\text{max}(j,k)]\in \B_n$,
which is already computed to be equal to
$$\relspinc(\x)+(j-k+n(i-j))\text{PD}[\mu].$$
Since $[\x,i,j,k,l]\in \D_1^s$ is mapped to both
$[\x,i,j]\otimes T^{-s}$
and to $[\x,l,k]\otimes T^{-s-n}$, it is
easy to check that the relative $\SpinC$ structure associated
with $[\x,i,j]\otimes T^k\in \B[n\ell]$ should be defined
equal to $s(\x)+(-k+n(i-j))\text{PD}[\mu]$. Here we assume
that $-\frac{n\ell}{2}\leq k<\frac{n\ell}{2}$. Note that for
any generator $\x\in \Ta\cap\Tb$, by definition we have
$s(\x)+i(\relspinc(\x))\text{PD}[\mu]=\relspinc(\x)$. It is
straight-forward to check that these assignments are respected
by the chain maps.\\

From the above constructions we get a map from
$\CFKT^\infty(Y_n(K),K_n,\relspinct)$ to the complex
$\M^{[s]}(\fbar_\ell)[\relspinct]$; the subset of
$\M^{[s]}(\fbar_\ell)$ corresponding to the relative
$\SpinC$ class $\relspinct$ with $[i(\relspinct)]=[s]$ in
$\frac{\Z}{n\Z}$. This map respects the filtration on the
two sides and it will produce a
quasi-isomorphism in homology in the following sense.\\

Recall that if $f:C\lra C'$
is a chain map between $\Z\oplus \Z$ chain complexes
such that for any positive test domain $P$ (see the
introduction for the definition of a positive test domain)
then $f$ will induce a chain map $f:C^P\lra C'^P$ where
$C^P$ and $C'^P$ denote the chain complexes generated by
the pre-images of $P$ in $C$ and $C'$. If
the induced map
$$f_*^P:H_*^P(C)\lra H_*^P(C')$$ is an isomorphism
for all positive test domains
we will call $f$ a \emph{quasi-isomorphism}. The above argument
proves that the constructed map
$$\CFK^\infty(Y_n(K),K_n,\relspinct)\lra
\M^{[s]}(\fbar_\ell)[\relspinct]$$
is a quasi isomorphism for any
$\relspinct\in \RelSpinC(Y_n(K),K_n)$, where $[s]=
[i(\relspinct)]\in \frac{\Z}{n\Z}$ is the class of
$i(\relspinct)$ modulo the integer $n$, and
$\M^{[s]}(\fbar_\ell)[\relspinct]$ is the subcomplex
of $\M^{[s]}(\fbar_\ell)$ produced by generators in
the relative $\SpinC$ class $\relspinct$. \\

In fact, for any positive test domain $P$ we may follow the
above process using complexes determined by the test domain
$P$ and the filtrations. Then we will obtain maps
$$(\CFKT^\infty(Y_n(K),K_n,\relspinct)^P\lra
(\M^{[s]}(\fbar_\ell)[\relspinct])^P$$
which give isomorphisms in homology.\\

As $m$ goes to infinity we may construct the stabilization of this
complex as follows. $\D_1^s(\relspinct)$ will be defined as before
without any restriction on how big $s$ is. The complex
$\D_1(\relspinct)$ will decompose as
$$\D_1(\relspinct)=\bigoplus_{\substack{s\in \Z\\ s=i(\relspinct)\
    (\text{mod }n)}} \D_1^s(\relspinct).$$
 We also introduce a complex $\B[\infty]=\B\otimes_\Z \Z[T,T^{-1}]$
as the complex generated by the elements
$\{[\x,i,j]\otimes T^k\ |\ \x\in \Ta\cap\Tb,\ i,j,k\in \Z\}$ over
$\Z$, with the differential
$$\partial_\B[\x,i,j]\otimes T^k=\sum_{\substack{\y\in \Ta\cap\Tb\\
\phi \in \pi_2(\x,\y)\\ \mu(\phi)=1}}\#(\widehat{\Mod}(\phi))
[\y,i-n_u(\phi),j-n_v(\phi)]\otimes T^k.$$
The maps $\fbar_\ell$ stabilize to give a map
$\fbar:\D_1\rightarrow \B[\infty]$ where $\fbar=f+v$ and
\begin{displaymath}
\begin{split}
&h[\x,i,j,k,l]=[\x,i,j]\otimes T^{-s},\ \ \ \ \ \ \ \\
&v[\x,i,j,k,l]=\tau[\x,l,k]\otimes T^{-s-n},\ \ \text{for }[\x,i,j,k,l]\in \D_1^s.
\end{split}
\end{displaymath}
 We may take the mapping cone $\M(\fbar)$ of $\fbar$.

The filtration $\Gi$ on $\M(\fbar)$ is
given by
\begin{displaymath}
\begin{split}
&\Gi:(\text{Generators of }\M(\fbar))\lra \Z\oplus\Z\\
&\begin{cases}
\Gi([\x,i,j,k,l])=(\text{max}(i,l),\text{max}(j,k))\ \ \
\ [\x,i,j,k,l]\in \D_1, \\
\Gi([\x,i,j]\otimes T^k)=(i,j)\ \ \ \
\ \ \  \ \ [\x,i,j]\otimes T^k\in \B\otimes_\Z\Z[\frac{\Z}{m\Z}]
\end{cases}
\end{split}
\end{displaymath}
and the decomposition of $\M(\fbar)$ according to the relative
$\SpinC$ classes  is determined by
assigning to
$a\in \M(\fbar)$ the relative $\SpinC$ class
$\relspinc(a)\in \RelSpinC(Y_n(K),K_n)$ defined by
\begin{displaymath}
\relspinc(a)=
\begin{cases}
\relspinc(\x)+((j-k)+n(i-j))\text{PD}[\mu]\ \ \ \ \ \ \ \text{if}\
a=[\x,i,j,k,l]\in \D_1\\
s(\x)+(-k+n(i-j))\text{PD}[\mu]\ \ \ \ \ \ \ \text{if}\
a=[\x,i,j]\otimes T^k \in \B[\infty].
\end{cases}
\end{displaymath}

The results of the above considerations, together with
the argument given in \cite{OS-Zsurgery} for theorem 4.1
will prove the following theorem:\\

\begin{thm}
\label{thm:Main} Suppose that $(Y,K)$ is a null-homologous knot,
$n$ is a nonzero integer, and that $(Y_n(K),K_n)$ is the result of
$n$-surgery on $(Y,K)$. Suppose that $\relspinct\in
\RelSpinC(Y_n(K),K_n)$ is a fixed relative $\SpinC$ class, and
that the complex $\M(\fbar)$ is
constructed as above. Then
$\CFKTT^\infty(Y_n(K),K_n,\relspinct)$ and $\M(\fbar)[\relspinct]$
are quasi-isomorphic.

\end{thm}
\begin{proof}
The proof is almost identical with the proof of theorem 4.1 in
\cite{OS-Zsurgery}, which occupies most of that paper. For the
reader's convenience, we sketch the proof omitting the
details.\\

The triangle maps of the previous section give rise to a triangle of
chain maps
\begin{displaymath}
\begin{split}
&f_1:\CFKT^\infty(Y_n(K),K_n)\lra \CFKT^\infty(Y_{m+n}(K),K_{m+n})\\
&f_2:\CFKT^\infty(Y_{m+n}(K),K_{m+n})\lra \bigoplus^m
  \CFKT^\infty(Y,K)\\
&f_3:\bigoplus^m \CFKT^\infty(Y,K) \lra \CFKT^\infty(Y_n(K),K_n)\\
\end{split}
\end{displaymath}
The maps $f_1$ and $f_2$ are constructed from the Heegaard diagram
$R_{m,n}$ as before. The map $f_3$ for a generator of the form
$[\x,i,j]T^k$ is defined by counting triangles $\psi
\in \pi_2(\x,\Theta,\y)$ (where $\Theta$ is the canonical generator
for the connected sum of $S^1\times S^2$s), such that the number
$m_p(\psi)$ is congruent to $c-k$ modulo the integer $m$, where
$c$ is a constant determined from the Heegaard diagram as in
\cite{OS-Zsurgery}. \\

Since the compositions $f_2\circ f_2$ and $f_3\circ f_2$ are
null-homotopic,
$f_1$ and $f_3$ induce maps $\Psi:\CFKT^\infty(Y_n(K),K_n)
\rightarrow \M(f_2)$
and $\Psi':\M(f_2) \rightarrow \CFKT^\infty(Y_n(K),K_n)$,
where $\M(f_2)$ denotes
the mapping cone of the filtered chain map $f_2$. These maps
are quasi-isomorphisms in the sense discussed earlier (i.e.
they induce isomorphisms in homology for any positive test
domain $P$).\\

Suppose that $m=\ell n$ for some large $\ell$.
Note that the chain map $f_2$ is  determined as the chain map
associated with the $4$-manifold cobordism
$$W_{\ell}:Y_{(\ell+1)n}(K)\lra Y.$$
For any $\SpinC$ structure $\spinct\in \SpinC(Y_{n(\ell+1)}(K))$
let $\xfrak_{\spinct},\yfrak_{\spinct}\in \SpinC(W_\ell)$
 be defined as before.
Then if we decompose $f_2$ as a sum of maps $$f_2=\sum_{\relspinct
\in \SpinC(W_\ell)}F_{W_\ell,\relspinct}\otimes T^{\m(\relspinct)},$$
and $a$ is a generator in $\SpinC$ class $\spinct\in
\SpinC(Y_{n(\ell+1)}(K)$ then
the only non-trivial components of $f_2(a)$ in this expression
(for large values of $\ell$) are the ones of the form
$F_{W_\ell,\xfrak_{\spinct}}(a)$ and $F_{W_\ell,\yfrak_{\spinct}}(a)$.
For any positive test domain $P$ one can check the commutativity of
the following diagrams:
\begin{displaymath}
\begin{array}{ccc}
\CFKT^P(Y_{n(\ell+1)}(K),K_{n(\ell+1)},\relspinct)
&\xrightarrow{F_{W_\ell,\xfrak_{\spinct}}}
& (\bigoplus^\ell \CFKT^\infty(Y,K))^P\\
\downarrow& &\downarrow\\
\D_1^s[\relspinct]^P&\xrightarrow{\ \ \ \ h\ \ \ \ }
&(\B[n\ell],\relspinct)^P
\end{array}
\end{displaymath}
and the diagram
\begin{displaymath}
\begin{array}{ccc}
\CFKT^P(Y_{n(\ell+1)}(K),K_{n(\ell+1)},\relspinct)
&\xrightarrow{F_{W_\ell,\yfrak_{\spinct}}}
& (\bigoplus^\ell \CFKT^\infty(Y,K))^P\\
\downarrow& &\downarrow\\
\D_1^s[\relspinct]^P&\xrightarrow{\ \ \ \ v\ \ \ \ }
&(\B[n\ell],\relspinct)^P.
\end{array}
\end{displaymath}
Here $\relspinct \in \RelSpinC(Y,K)$ is a relative $\SpinC$
structure, $\spinct$ is the induced $\SpinC$ structure on
$Y_{n(\ell+1)}(K)$ and $-\frac{n(\ell+1)}{2}\leq s
<\frac{n(\ell+1)}{2}$ is an integer with  the property that
$i(t)=s$ modulo the integer $n(\ell+1)$. Furthermore,
$(\B[m],\relspinct)$ is the part of $\B[m]$ in the relative
$\SpinC$ class $\relspinct$.\\

Note that for a fixed relative $\SpinC$ class $\relspinct\in
\RelSpinC(Y,K)$
$$\fbar_\ell :\bigoplus_{\substack{
-\frac{n(\ell+1)}{2}\leq s< \frac{n(\ell+1)}{2}\\ s=i(\relspinct)
(\text{ mod }n)}} \D_1^s(\relspinct)\lra \bigoplus^\ell
(\B[m],\relspinct)$$ will produce a mapping cone $\M^{[s]}(\fbar_\ell)
[\relspinct]$ which is quasi-isomorphic to  the mapping cone of $f_2$ in the
relative $\SpinC$ class $\relspinct$. This may be verified directly
from the above commutative diagrams.\\

For $s\geq \frac{n(\ell+1)}{2}$ note that $i(\relspinc(\x))+j-k=s$
implies that for the generators $[\x,i,j,k,l]$ in $\D_1^s$ we have
$j=k+(s-i(\relspinc(\x)))$ and $i=l+(s-i(\relspinc(\x)))+1$. This
implies that if $\ell$ is large enough  $i$ is the maximum of
$\{i,l\}$ and $j$ is the maximum of $\{j,k\}$. This implies that
the map $h$ from $(\D_1^s)^P$ to its image in $\B[\infty]$ is an
isomorphism. For such values of $s$, one may consider $\phi=v\circ
h^{-1}:\D_1^s\lra \D_1^{s+n}$ and observe that $\phi^p([\x,i,j,k,l]$
will have a $\Z^4$-filtration (coming from the integer components)
which is less than or equal to
$$\Big[pl-pi+i-\Big(\begin{array}{c}p\\ 2\end{array}\Big),
pk-pj+j-\Big(\begin{array}{c}p\\ 2\end{array}\Big),
pk-pj+k-\Big(\begin{array}{c}p+1\\ 2\end{array}\Big),
pl-pi+i-\Big(\begin{array}{c}p+1\\ 2\end{array}\Big)
\Big].$$
The positivity of $P$ implies that for fixed $(i,j,k,l)\in \Z^4$
this sequence will eventually leave $P$.

If $s\leq -\frac{n(\ell+1)}{2}$ a similar argument
shows that $v$ from $(\D_1^s)^P$ to its image in $\B[\infty]$ is an
isomorphism and $(h\circ v^{-1})^p(a)$ will eventually vanish
for a fixed generator $a$ and large enough $p$.
These observations imply that for large values of
$\ell$, and any positive test domain $P$ there is an isomorphism
induced by the inclusion
$$H_*^P(\M^{[s]}(\fbar_\ell)[\relspinct])\lra H_*^P(\M(\fbar)[\relspinct]),$$
where $[s]$ denote the class of $i(\relspinct)$ modulo the integer
$n$. Thus the chain complex $\M(\fbar)[\relspinct]$ is
quasi-isomorphic to $\CFKT^\infty(Y_n(K),K_n,\relspinct)$.
This would complete the proof.
\end{proof}

Suppose that $K$ is a knot in a homology sphere $Y$ and $s\in
\Z=\RelSpinC(Y,K)$. Denote by $\A_s^+$ the complex generated by
those generators $[\x,i,j]\in (\Ta\cap\Tb)\times \Z^2$ such that
max$(i,j)\geq 0$, and $i(\x)+i-j=s$ as discussed in
\cite{OS-Zsurgery}, where $i(\x)=i(\relspinc(\x))$. Let
${\bar{\A}}^+_s$ denote the complex generated by generators
$[\x,i]\in (\Ta\cap\Tb) \times \Z^{\geq 0}$ giving $\CFT^+(Y)$ for
each $s\in \Z$. Let $\mathcal{A}_s^+$ denote the direct sum
$\mathcal{A}^+_s= \bigoplus_{t= s\ (\text{mod}\ n)}\A^+_t$ and
similarly define $\bar{\mathcal{A}}^+_s$. Let
$h',v':\mathcal{A}^+_s \lra \bar{\mathcal{A}}^+_s$ be the maps
sending $[\x,i,j]\in {\A_t}^+$ to
\begin{displaymath}
\begin{split}
&h'[\x,i,j]=[\x,i]\in{\bar{\A}}^+_t,\ \ \text{and}\\
&v'[\x,i,j]=\tau[\x,j]\in {\bar\A}^+_{t-n}
\end{split}
\end{displaymath}
respectively. Denote by $\zeta_s$ the sum of these two maps.
The above theorem implies the result of \cite{OS-Zsurgery}
that the chain homotopy type of $\CFT^+(Y_n(K),s)$ is the
same as that of the mapping cone of $\zeta_s$:\\

\begin{cor}\label{cor:KtoY}
Suppose that $(Y,K)$ and $Y_n(K)$ are the same as before. Define
$\mathcal{A}^+_s$ and $\bar{\mathcal{A}}_s^+$ and also the map
$\zeta_s$ as above for $-\frac{n}{2}\leq s< \frac{n}{2}$. Then the
homology groups $\text{\emph{HF}}^+(Y_n(K),\spinc_n)$ (where
$\spinc_n$ is any $\SpinC$ structure with $i_n(\spinc_n)=s
(\text{mod}\ n)$) is the same as the homology of the mapping cone
of $\zeta_s$ in the $\SpinC$ class $\spinc_n$, as computed in
\cite{OS-Zsurgery}.
\end{cor}
\begin{proof}
If the three-manifold $Y$ is a homology sphere and $(Y,K)$ is
the given knot, then the set of relative $\SpinC$ structures
$\RelSpinC(Y_n(K),K_n)=\RelSpinC(Y,K)$ is naturally isomorphic
with $\Z$. The map
$$\Gi_n:\Z=\RelSpinC(Y_n(K),K_n)\lra \SpinC(Y_n(K))=\frac{\Z}{n\Z}$$
is the reduction modulo the integer $n$. For
$-\frac{n}{2}\leq s<\frac{n}{2}$, thought of as an element of
$\RelSpinC(Y_n(K),K_n)$, the complex
$\CFKT^\infty(Y_n(K),K_n;s)$ is a $\Z\oplus\Z$-filtered chain complex
which is, as a $\Z$ filtered chain complex, chain homotopic with
$\CFT^\infty(Y_n(K),[s]_{\text{mod }n})$. Using theorem~\ref{thm:Main}
$C^+(s)=\CFT^+(Y_n(K),[s]_{\text{mod }n})$ is generated by two types of
generators, some of them in $\D_1$ and some of them in $\B[\infty]$.\\

If $[\x,i,j,k,l]\in \D_1$ is a generator of  $C^+(s)$
it is implied that $i(\x)+(j-k)+n(i-j)=s$. This implies
that if $i(\x)+j-k$ is congruent to $s$ modulo $n$, then
$i=-\frac{i(\x)+j-k-s}{n}+j$ will compute $i$. Note that
from the relation $i-j+k-l=1$, the integer $l$ may be computed as well.
The filtration $\Gi$ gives the value
$(\text{max}(i,l),\text{max}(j,k))$ on this generator.
As a result the projection over the second component of the
$\Z\oplus\Z$ filtration will be
$\text{max}(j,k)$. Let $\A_t^+$ for $t= s\ (\text{mod}\ n)$
denote the complex generated by those $[\x,j,k]\in (\Ta\cap\Tb)\times \Z^2$
such that $i(\x)+j-k=t$ and max$(j,k)\geq 0$. There is a natural
correspondence between the generators of $C^+(s)$ of the form
$[\x,i,j,k,l]\in \D_1$ with the generators of $\bigoplus_{t= s\
  (\text{mod}\ n)}\A^+_t$.\\

If $[\x,i,j]\otimes T^k \in \B[\infty]$ is a generator of $C^+(s)$
it is implied that $-k+n(i-j)=s$. As a result, $-k= s\
(\text{mod}\ n)$, and the value of $i$ is determined once such a
choice for $k$ is fixed. Thus, the generators of the form
$[\x,i,j]\otimes T^{-t}$ for $t= s\ (\text{mod}\ n)$ are in
correspondence with generators of $\bar{\A}^+_t$. This complex is
generated by $[\x,i]\in (\Ta\cap\Tb)\times \Z$ with $i\geq 0$,
which gives $\CFT^+(Y)$. $\bar\A^+_t$ is the same complex
for all values of $t$.\\

Under the above correspondence the map $h$ takes $[\x,j,k]\in
\A^+_t$ to $[\x,j]\in \bar\A^+_t$, and the map $v$ takes the same
generator to $\tau [\x,k]\in \bar\A^+_{t-n}$, where $\tau$ is the
chain homotopy equivalence discussed earlier. The homology of the
mapping cone $C^+(s)$ is thus computed
exactly as proved in \cite{OS-Zsurgery}.\\
\end{proof}

Finally note that in the case of null-homologous knots, there is a
natural isomorphism of chain complexes (without the filtration)
from $\D_1$ to $\B[\infty]$. Namely, we may define:
\begin{displaymath}
\begin{split}
&\Psi:\D_1 \lra \B[\infty]\\
&\Psi[\x,i,j,k,l]:=[\x,i,j]\otimes T^{-i(\relspinc(\x))+k-j}.\\
\end{split}
\end{displaymath}
It is easy to verify that this in fact is an isomorphism. The map
$\fbar$ will introduce a map $\gbar_n:\D^{up}\lra \D^{down}$, where
$\D^{up}$ and $\D^{down}$ are copies of $\D_1$. The map induced from
$h$ is the identity. The map $v_n$ induced by $v$ is more interesting:
$$v_n[\x,i,j,k,l]=\tau[\x,l,k,2k-j-n,2l-i-n].$$
Note that the chain homotopy equivalence on $\B$ may be naturally
extended to a chain homotopy equivalence from the complex $\D_1$
(equipped with the $\Z\oplus\Z$ filtration coming from projection
over the $3$rd and $4$th integer components of the generators) to
$\D^{down}$. There is a discussion on this in the introduction. We
denote this new chain homotopy equivalence by the same letter $\tau$.\\

The filtration induced on $\D^{down}$ from $\B[\infty]$ is given by
the first two integer components of the generators. We obtain the
following re-statement of theorem~\ref{thm:Main}:

\begin{thm}\label{thm:Re-Main}
Let $(Y,K),n$ and $\D_\delta$ be as above. Let $\D^{up}$ and
$\D^{down}$ be copies of the complex $\D_1$. Denote by
$v_n:\D^{up}\lra \D^{down}$ the map defined by $v_n[\x,i,j,k,l]=\tau
[\x,l,k,2k-j-n,2l-i-n]$. Let $\gbar_n=Id+v_n$ and denote by
$\M(\gbar_n)$ the mapping cone of $\gbar_n$. Define a filtration
$\Gi$ on $\M(\gbar_n)$ by setting
\begin{displaymath}
\Gi[\x,i,j,k,l]=
\begin{cases}(\text{max}(i,l),\text{max}(j,k))
\ \ \ \text{if}\ [\x,i,j,k,l]\in \D^{up},\\
(i,j)\ \ \ \ \ \ \ \ \ \ \ \ \ \ \ \ \ \ \
\ \
\ \text{if } [\x,i,j,k,l]\in \D^{down}.\\
\end{cases}
\end{displaymath}
Also define a map from the set of generators to the set of relative
$\SpinC$ structures $\RelSpinC(Y_n(K),K_n)=\RelSpinC(Y,K)$ by
\begin{displaymath}
\begin{split}
\relspinc[\x,i,j,k,l]=
\relspinc(\x)+((j&-k)+n(i-j))\text{\emph{PD}}[\mu]\\
&\text{for any generator } [\x,i,j,k,l]\in \D^{up}\ \text{or}\ \D^{down},\\
\end{split}
\end{displaymath}
and split $\M(\gbar_n)$ as $\M(\gbar_n)
=\bigoplus_{\relspinct\in
\RelSpinC(Y,K)}\M(\gbar_n)[\relspinct]$.
Then $\CFKTT^\infty(Y_n(K),K_n,\relspinct)$
is quasi-isomorphic to
$\M(\gbar_n)[\relspinct]$ for every
$\relspinct\in \RelSpinC(Y,K)$.
\end{thm}
It will be more convenient to state our next result as a
generalization of the first form of this theorem. The
theorem look nicer, however, in this second form when we
work with a null-homologous knot.\\

\section{Rationally null-homologous knots}
In this section we generalize the construction of previous section to
the case of rationally null-homologous knots $(Y,K)$. We remind the
reader of a couple of facts from \cite{OS-Qsurgery} where the notion
of knot Floer homology is generalized to the case of rationally
null-homologous knots and also the integral surgery
formulas are generalized to Morse surgery formulas for this type of
knots.\\

Note that to a rationally null-homologous knot $(Y,K)$ is associated
a notion of relative $\SpinC$ structure and the set of such structures
is denoted by  $\RelSpinC(Y,K)$. There is a surjective reduction map
$$G_{Y,K}:\RelSpinC(Y,K)\lra \SpinC(Y).$$
If $H=(\Sig,\alphas,\betas,p)$ is a Heegaard diagram for $(Y,K)$ as
before, there is a map associated with the marked point $p$ which
assigns relative $\SpinC$ structures to intersections of $\Ta$ and
$\Tb$. Namely we have the map
$$\relspinc:\Ta\cap\Tb\lra \RelSpinC(Y,K).$$
The choice of a framing $\lambda$ determines a push-off $K_\lambda$
of the knot $K$ into the knot complement $Y\setminus \text{nd}(K)$,
which provides us, via Poincar\'e duality, with a cohomology class
$$\text{PD}[\lambda]:=\text{PD}[K_\lambda]\in
H^2(Y\setminus\text{nd}(K),\partial(Y\setminus \text{nd}(K)),\Z).$$
Note that the set of relative $\SpinC$ structures is an affine space
over this later cohomology group.\\

 If $\lambda$ is a framing
for $(Y,K)$, then $\lambda+n\mu$ is also a framing, and we may define
the push-off $K_{\lambda+n\mu}$ similar to $K_\lambda$. We may also
define the cohomology class
PD$[\lambda+n\mu]=$PD$[\lambda]+n$PD$[\mu]$ similarly.
\\

We may use the framing in place of the curve $\lambda$ in previous
sections to define the curves $\lambda_n$, and the knots
$(Y_n(K),K_n)$. Again
the set of relative $\SpinC$ structures associated with a knot
$(Y,K)$ is the same as the set of relative $\SpinC$ structures
associated with $(Y_n(K),K_n)$:
$$\RelSpinC(Y,K)\cong \RelSpinC(Y_n(K),K_n).$$

Fixing the framing $\lambda$ we may start the process
of second and third sections.  Suppose that a
Heegaard diagram $H=(\Sig,\alphas,\betas,p)$ for $(Y,K)$ is given
as above, inducing a differential $\partial^\infty$ on
$\CFKT^\infty(Y,K)$, which is generated by $(\Ta\cap\Tb)\times \Z\times
\Z$.
Let $\D$ be the complex generated by $(\Ta\cap\Tb)\times \Z^4$ with
the differential
$$\partial_\D[\x,i,j,k,l]=\sum_pn_p[\y_p,i-i_p,j-i_p,k-k_p,l-k_p],$$
where $\partial^\infty[\x,i,k]=\sum_pn_p[\y_p,i-i_p,k-k_p]$. For
$\delta\in \Z$ let $\D_\delta$ denote the subcomplex of $\D$
generated by the generators $a=[\x,i,j,k,l]$ such that
$\Delta(a)=i-j+k-l=\delta$.\\

  In
\cite{OS-Qsurgery} a map $\Xi:\SpinC(Y_n(K))\lra \RelSpinC(Y_n(K),K_n)$
is constructed (in the presence of the framing $\lambda$ and for large
values of $n$) which plays the role of the map sending
$\spinct\in \SpinC(Y_n(K))$ to $\relspinct\in \RelSpinC(Y_n(K),K_n)$
such that $s_n(\relspinct)=\spinct$ and $-\frac{n}{2}\leq
i(\relspinct)< \frac{n}{2}$. \\

Similar to the definition of $\D_1^s(\relspinct)$ for any relative
$\SpinC$ class $\relspinct\in \RelSpinC(Y_n(K),K_n)$ and any
$\SpinC$ class $\spinct_n\in \SpinC(Y_n(K))$ denote by
$\D_1^{\spinct_n}(\relspinct)$ the subcomplex of $\D_1$ generated by
the generators $[\x,i,j,k,l]\in \D_1$ satisfying
\begin{displaymath}
\begin{cases}
\relspinc(\x)+(j-k){\text{{PD}}}[\mu]+(i-j)\text{{PD}}
[\lambda+n\mu]=\relspinct\\
\relspinc(\x)+(j-k)\text{PD}[\mu]=\Xi(\spinct_n).
\end{cases}
\end{displaymath}
It is implied that if $\D_1^{\spinct_n}(\relspinct)$ is non-empty
then $G_n(\relspinct)=\spinct_n$.\\

The complexes $\D_1^{\spinct_n}(\relspinct)$ are the natural
replacements for $\D_1^s(\relspinct)$, and we may follow the
process used for proving theorem~\ref{thm:Large-Surgery} to prove
the following:

\begin{thm}
\label{thm:QLarge-Surgery}
Suppose that $(Y,K)$ is a rationally null-homologous knot, $\lambda$
is a framing for $K$ and $(Y_n(K),K_n)$ is as above. Construct
the complex $\D_1$ as before
Then for large values of $n\in \Z$ the
filtered chain complex associated with the rationally
null-homologous knot $(Y_n(K),K_n)$ in relative $\SpinC$ class
$\relspinct\in \RelSpinC(Y_n(K),K_n)$ has the same chain homotopy
type as the complex $\D_1^{\spinct_n}(\relspinct)$
equipped with the $\Z\oplus\Z$
filtration given by
$$\mathcal{F}[\x,i,j,k,l]=(\text{max}\{i,l\},\
\text{max}\{j,k\}).$$
Here we have chosen $\spinct_n$ so that $\spinct_n=G_n(\relspinct)$.
\end{thm}
\begin{proof}
All the steps in the proof are completely similar to the steps
in the proof for the null-homologous case.
\end{proof}

Define the complex $\B[\infty]$ as the complex generated by the
generators $[\x,i,j]\otimes T^{\relspinct}$ where $\x\in
\Ta\cap\Tb$, $i,j\in \Z$ and $\relspinct\in \RelSpinC(Y,K)$ and we
have the relation
$$\spinc(\x)=G_{Y,K}(\relspinct)=G(\relspinct)$$
Clearly this is a generalization of the definition of the complex
$\B[\infty]$ used in the third section.\\

We  may construct two maps from the complex
$\D_1$ to $\B[\infty]$ as follows. These maps will be given via
the formulas
\begin{displaymath}
\begin{split}
&h,v_\lambda:\D_1\lra \B[\infty],\\
&h[\x,i,j,k,l]=[\x,i,j]\otimes T^{\relspinc(\x)
+(j-k)\text{PD}[\mu]}\\
&v_\lambda[\x,i,j,k,l]=\tau[\x,l,k]\otimes T^{\relspinc(\x)
+(j-k)\text{PD}[\mu]+\text{PD}[\lambda]}.\\
\end{split}
\end{displaymath}
Define $f_\lambda=h+v_\lambda$ and let $\M(f_\lambda)$ denote the
mapping cone of $f_\lambda$. Define a filtration $\Gi$ on the
generators of $\M(f_\lambda)$ by
\begin{displaymath}
\begin{cases}
\Gi([\x,i,j,k,l])=(\text{max}(i,l),\text{max}(j,k))\ \ \ [\x,i,j,k]\in
\D_1\\
\Gi([\x,i,j]\otimes T^{\relspinct})=(i,j)\ \ \ \ \ \ \
 \ \ \ \ \ \ \ \
   [\x,i,j]\otimes T^{\relspinct}\in \B[\infty].
\end{cases}
\end{displaymath}
The relative $\SpinC$ classes of generators in $\M(f_\lambda)$ will
be defined via
\begin{displaymath}
\relspinc(a)=
\begin{cases}
\relspinc(\x)+(j-k)\text{PD}[\mu]+(i-j)\text{PD}[\lambda]\ \ \text{ if
} a=[\x,i,j,k,l]\in \D_1,\\
\relspinct+(i-j)\text{PD}[\lambda]\ \ \ \ \ \ \ \ \ \
\ \ \ \ \ \ \ \ \ \ \ \
\text{ if }a=[\x,i,j]\otimes T^{\relspinct} \in \B[\infty].
\end{cases}
\end{displaymath}

 We may insert these
constructions, which are the generalized versions of
the previous ones, in the proof of theorem~\ref{thm:Main}
to obtain the following.

\begin{thm}
\label{thm:QMain} Suppose that $(Y,K)$ is a rationally
null-homologous knot, $\lambda$ is a framing for $K$, and that
$(Y_\lambda(K),K_\lambda)$ is the knot obtained as above by Morse
surgery with framing $\lambda$ on $(Y,K)$. Suppose that the
complex $\M(f_\lambda)$ be as above, and let $\relspinct\in
\RelSpinC(Y_n(K),K_n)$ be a relative $\SpinC$ structure. Then the
$\Z\oplus \Z$-filtered chain complex
$\CFKTT^\infty(Y_\lambda(K),K_\lambda,\relspinct)$  is
quasi-isomorphic to the complex $\M(f_\lambda)[\relspinct]$
associated with the mapping cone of the map
$f_\lambda=h+v_\lambda:\D_1\lra \B[\infty]$, equipped with the
$\Z\oplus \Z$-filtration $\Gi$.
\end{thm}
\begin{proof}
Again, all the necessary modifications are minor. The proof of
theorem~6.1 in \cite{OS-Qsurgery} may be combined with our techniques
for the proof of theorem~\ref{thm:Main} to prove the above theorem.
\end{proof}
\begin{remark}
This is a generalization of theorem 6.1 from \cite{OS-Qsurgery} in
an obvious way.
\end{remark}

\section{Rational surgeries on null-homologous knots}
The natural application of the above generalization of
theorem~\ref{thm:Main} is a computation for null-homologous knots
of the homologies of the Heegaard Floer complex associated
with rational surgeries on them. \\

For simplicity, we choose to deal with the case where $Y$ is a
homology sphere, so that the knot $(Y,K)$ is automatically
null-homologous.
As in \cite{OS-Qsurgery}, if $\frac{p}{q}\in \mathbb{Q}$
is a rational number,
then write
$$\frac{p}{q}=\frac{r}{q}+a=\frac{r}{q}+\lfloor\frac{p}{q}\rfloor,$$
and note that $(Y_{\frac{p}{q}}(K),K_{\frac{p}{q}})$ may
be obtained by a Morse surgery with coefficient $a$ on the
knot $K\#O_{\frac{q}{r}}\subset Y\#L(q,r)$. We remind the reader
that $O=O_{\frac{q}{r}}$ is the knot obtained as one component of
the Hopf link in the three-manifold $L=L(q,r)$ obtained by a
$\frac{q}{r}$ surgery on the second component of the Hopf link.
Note that $(L,O)$ is a $U$-knot, according to \cite{OS-Qsurgery}.\\

We remind the reader of a number of facts from \cite{OS-Qsurgery}
about the splitting of relative $\SpinC$ structures under connected
sum of knots and about the filtered chain homotopy type of
$(Y_1\#Y_2,K_1\#K_2)$ (in terms of the chain homotopy type of
$(Y_1,K_1)$) when $(Y_2,K_2)$ is a $U$-knot. First note that there
is a connected sum map
$$\RelSpinC(Y_1,K_1)\times \RelSpinC(Y_2,K_2)\lra
\RelSpinC(Y_1\#Y_2,K_1\#K_2)$$
sending  a pair of relative $\SpinC$ structures
$(\relspinc_1,\relspinc_2)$ to $\relspinc_1\#\relspinc_2$.
Note that there is a one-parameter family of pairs $(\relspinc_1,\relspinc_2)$
such that $\relspinc_1\#\relspinc_2$ is a fixed relative
$\SpinC$ class in $\RelSpinC(Y_1\#Y_2,K_1\#K_2)$. If $K_2$ is a
$U$-knot for any $\relspinc_1\in \RelSpinC(Y_1,K_1)$ and any
$\spinc_2\in \SpinC(Y_2)$ there exists a unique relative
$\SpinC$ class $\relspinc_2\in \RelSpinC(Y_2,K_2)$ with the
property that $G_{Y_2,K_2}(\relspinc_2)=\spinc_2$ such that there is
an equivalence of chain homotopy types
$$\CFKT^\infty(Y_1,K_1,\relspinc_1)\cong \CFKT^\infty
(Y_1\#Y_2,K_1\#K_2,\relspinc_1\#\relspinc_2).$$
In  particular for the knot $(L(q,r),O_{\frac{q}{r}})$
there is a commutative diagram
\begin{displaymath}
\begin{array}{ccc}
\Z&\xrightarrow{\ \ \ \phi\ \ \ }&\RelSpinC(L(q,r),O_{\frac{q}{r}})\\
\downarrow& &\downarrow \text{}_{\text{}^{G_{L,O}}}\\
\frac{\Z}{q\Z}&\xrightarrow{\ \ \ \cong\ \ \ } &\SpinC(L(q,r))
\end{array}
\end{displaymath}
such that for $0\leq i \leq q-1$ there is an isomorphism
$\widehat{\HFKT}(L(q,r),O_{\frac{q}{r}};\phi(i))\cong \Z$, and
such that for all other $i$ we have
$\widehat{\HFKT}(L(q,r),O_{\frac{q}{r}},\phi(i))=0$.\\

For any homology sphere $Y$ and any knot $(Y,K)$, we may note that
\begin{displaymath}
\begin{split}
&\RelSpinC(Y,K)\cong \RelSpinC(Y\#L,K\#O)\cong \Z\ \text{and}\\
&H^2(Y,K)\cong H^2(L,O)\cong H^2(Y\#L,K\#O)\cong \Z.
\end{split}
\end{displaymath}
Under these isomorphisms the following diagram is commutative
\begin{displaymath}
\begin{array}{ccc}
\Z\oplus \Z &\xrightarrow{\ \ \ f\ \ \ }&\Z\\
\downarrow& &\downarrow\\
H^2(Y,K)\oplus H^2(L,O)&\xrightarrow{\ \ \ \ \ } & H^2(Y\#L,K\#O)
\end{array}
\end{displaymath}
where $f$ is defined via $f(x,y)=qx+y$. Suppose that $K_\lambda$
is the push-off of the knot $K\#O$ with respect to the framing
$a$ (where $a+\frac{r}{q}=\frac{p}{q}$) into the complement of this
knot in $Y\#L$. Then according to \cite{OS-Qsurgery} the Poincar\'e
dual PD$[\lambda]$ of the homology class represented by $K_\lambda$
represents the element
$$p\in \Z=H^2(Y\#L,K\#O)=H^2(Y_{\frac{p}{q}}(K),K_{\frac{p}{q}}).$$
The meridian of the knot $K\#O$ in $Y\#L$ is just the image of the
meridian $\mu$ of the knot $(Y,K)$ in the connected sum. As a result,
the push-off $K_\mu$ is obtained as the image of the push-off of the curve
$\mu$ (in the complement of $K$ in $Y$) under the map $f$ constructed
above. Using the above isomorphisms, this corresponds to the element
$$q\in \Z=H^2(Y\#L,K\#O)=H^2(Y_{\frac{p}{q}}(K),K_{\frac{p}{q}}).$$

We need to construct the complex $\D_1$ out of the complex $\D$
associated with
$\CFKT^\infty(Y\#L,K\#O)$. To this end note that this complex
in the relative $\SpinC$ class $s\in \Z=\RelSpinC(Y\#L,K\#O)$
corresponds to the complex $\CFKT^\infty(Y,K)$ in the relative
$\SpinC$ class $[\frac{s}{q}]\in \Z=\RelSpinC(Y,K)$. More precisely
there is a filtered quasi-isomorphism between the two complexes. \\

Let $H=(\Sig,\alphas,\betas,p)$ be a pointed Heegaard diagram for
the knot $(Y,K)$ and assume that the complex $\CFKT^\infty(Y,K)$ is
generated by the generators $[\x,i,j]\in
(\Ta\cap\Tb)\times\Z\times\Z$. According to the above
paragraph, the complex $\CFKT^\infty(Y\#L,K\#O)$ is generated by
the generators of the form $[\x,i,j]\otimes \zeta^t$ where $0\leq
t< q$, $\x\in \Ta\cap\Tb$ and $i,j\in \Z$. This implies that the
corresponding complex $\D'$ associated with
$\CFKT^\infty(Y\#L,K\#O)$ is of the form $\D\otimes_\Z
\frac{\Z[\zeta]}{\zeta^q=1}$, where $\D$ is associated with the
Heegaard diagram $H$ as before. Similarly $\D'_1=\D_1\otimes_\Z
\frac{\Z[\zeta]}{\zeta^q=1}$ may be obtained. The $\Z\oplus\Z$
filtration on $\D'_1$ comes from the $\Z\oplus \Z$ filtration on
the $\D_1$ factor according to the above construction. The relative
$\SpinC$ structures associated with the generators of $\D_1'$ are
described as follows. If $\x\in \Ta\cap\Tb$ is a generator and if
$i(\x)\in \Z=\RelSpinC(Y,K)$ is the relative $\SpinC$ structure
associated with it then the relative $\SpinC$ structure in
$\Z=\RelSpinC(Y\#L,K\#O)=\RelSpinC
(Y_{\frac{p}{q}}(K),K_{\frac{p}{q}})$ associated with $\x\otimes
\zeta^t$ under the above correspondence will be $qi(\x)+t$. The
computation of relative $\SpinC$ structures in
theorem~\ref{thm:QMain} then implies that the relative $\SpinC$
class associated with a generator $[\x,i,j,k,l]\otimes \zeta^t$ is
given by the following formula
$$i([\x,i,j,k,l]\otimes \zeta^t)=qi(\x)+p(i-j)+q(j-k)+t\in \Z\cong
\RelSpinC(Y_{\frac{p}{q}}(K),K_{\frac{p}{q}}).$$

Note that $\B[\infty]$ is generated by generators $[\x,i,j]\otimes
\zeta^t \otimes T^{\relspinct}$ such that $qi(\x)+t=\relspinct$.
As a result, the value of $t$ is determined from $\relspinct\in
\Z\cong\RelSpinC(Y_{\frac{p}{q}},K_{\frac{p}{q}})$. Thus the
complex $\B^\infty$ is in fact generated by the generators of the
form $[\x,i,j]\otimes T^{\relspinct}$ where $\x \in \Ta\cap\Tb$,
$i,j\in \Z$ and $\relspinct\in \Z\cong
\RelSpinC(Y_{\frac{p}{q}},K_{\frac{p}{q}})$. The relative $\SpinC$
structure in $\RelSpinC(Y_{\frac{p}{q}},K_{\frac{p}{q}})$
associated with any such generator is given, according to
theorem~\ref{thm:QMain}, by the following formula
$$i([\x,i,j]\otimes T^{\relspinct})=\relspinct+p(i-j)\in \Z\cong
\RelSpinC(Y_{\frac{p}{q}},K_{\frac{p}{q}}).$$
The maps from $\D_1\otimes_\Z \frac{\Z[\zeta]}{\zeta^q=1}$ to
$\B[\infty]$ are given by
\begin{displaymath}
\begin{split}
&h([\x,i,j,k,l]\otimes\zeta^t)=[\x,i,j]\otimes T^{q(i(\x)+j-k)+t}\\
&v([\x,i,j,k,l]\otimes\zeta^t)=[\x,l,k]\otimes T^{q(i(\x)+j-k)+t+p}.
\end{split}
\end{displaymath}
Let $\M(\fbar)$ denote the mapping cone of $\fbar=h+v$.
Define a $\Z\oplus\Z$
grading on the generators of $\M(\fbar)$ by
\begin{displaymath}
\begin{split}
&\Gi([\x,i,j,k,l]\otimes\zeta^t)=
(\text{\emph{max}}(i,l),\text{\emph{max}}(j,k))\\
&\Gi([\x,i,j]\otimes T^{\relspinct})=(i,j).
\end{split}
\end{displaymath}
 This
complex (and consequently its homology) is decomposed into a direct
sum according to
the relative $\SpinC$ structures: $$\M(\fbar)=\bigoplus_{\relspinct
\in \RelSpinC(Y,K)} \M(\fbar)[\relspinct].$$

Although the filtered chain homotopy type may change in the course of this
process (as a quasi-isomorphism is composed with a chain homotopy
equivalence), the homology is preserved.

\begin{thm}\label{thm:QSF}
Let $Y$ be a homology sphere and let $(Y,K)$ denote a knot in
$K$. Suppose that $\frac{p}{q}> 0$ is a rational number and let
$(Y_{\frac{p}{q}},K_{\frac{p}{q}})$, $\fbar$ and $\M(\fbar)$ be as
before.
Then for any relative $\SpinC$ class $\relspinct\in \RelSpinC(Y,K)$,
the $\Z\oplus\Z$-filtered chain complex
$\CFKTT^\infty(Y_{\frac{p}{q}}(K),K_{\frac{p}{q}},\relspinct)$ is
quasi-isomorphic to the mapping cone $\M(\fbar)[\relspinct]$.
\end{thm}

Note that this
theorem may be re-stated as in the introduction.\\

\section{Non-vanishing results for $\widehat{\HFKT}(K_{\frac{p}{q}})$}

In this section we consider the special case where $P=\{(0,0)\}$
and $K$ is a knot in $S^3$ where the construction simplifies
significantly. A non-vanishing result may be proved for rational
surgeries on $K$ which may be used for re-proving \emph{Property
P}
as discussed in the introduction.\\

We will use the rational surgery formula as stated in the introduction.
Suppose that $\frac{p}{q}$ is a positive rational number and let
$L=K_{\frac{p}{q}}$ denote the result of $\frac{p}{q}$ surgery on
$K$. For simplicity, we will denote $\D^{up}$ by $\A$ and
$\D^{down}$ by $\B$. For this particular choice of $P$, we will
denote $\A^P$ by $\Ahat$ and $\B^P$ by $\Bhat$. The map
$f=I+g_{\frac{p}{q}}$ from $\A$ to $\B$ induces a chain map from
$\Ahat$ to $\Bhat$ and $\widehat{\HFKT}(L)$ is in fact the homology
of the complex
$$H_*(\Ahat)\xrightarrow{\ \ f_*\ \ }H_*(\Bhat).$$

Note that
$\Bhat(\relspinct)$ for any relative $\SpinC$ structure
$$\relspinct\in \Z\cong \RelSpinC(S^3,K)$$
is generated by $[\x,0,0,k,k-1]\otimes T^t$ such that
$q(i(\x)-k)+t=\relspinct$. This implies that
$i(\x)-k=\relspinc=\lfloor{\relspinct}/{q}\rfloor$ and
$t=q\{{\relspinct}/{q}\}$. The homology of this complex
is just $\widehat{\HFT}(S^3,\spinc_0)=\Z$, where $\spinc_0$ is
the unique $\SpinC$ structure on $S^3$.\\

The complex $\Ahat$ is more interesting. Consider a generator
$a=[\x,i,j,k,l]\otimes T^t$ in $\Ahat(\relspinct)$.
Suppose that $i=j+\delta$. It is implied that $l=k-1+\delta$. As a
result
$$\text{max}(j,k)=0=\text{max}(i,l)=\delta+\text{max}(j,k-1).$$
This can happen only if $\delta \in \{0,1\}$. Correspondingly we
may write $\Ahat=\Ahat_0\oplus \Ahat_1$, where $\Ahat_i$ is the
part of $\Ahat$ generated by generators as above such that $\delta=i$.\\

If $a$ is in $\Ahat_0$ then $i=j=0$ and $k\leq 0$. As a result the
complex $\Ahat_0$ may be identified by the complex $\Bhat\{k\leq 0\}$
consisting of the part of complex $\Bhat$ with non-positive
$k$-component. In fact $\Ahat_0(\relspinct)=\Bhat\{k\leq
0\}(\relspinct)$. If $C=\CFKT^\infty(K)$ denotes the
complex generated by $[\x,i,k]$
then one can check that in fact
$$\Ahat_0(\relspinct)\cong \Bhat\{k\leq 0\}(\relspinct)\cong
C\{i=0,k\leq 0\}(\lfloor\frac{\relspinct}{q}\rfloor).$$
The map $f_*$ will be the map induced in homology by the inclusion of
$\Ahat_0$ in $\Bhat$.\\

However, if $a$ is in $\Ahat_1$, then $k=l=0$ and $i\leq 0$. We will
also have
$$q(i(\x)+i-1)+p+t=\relspinct.$$
Thus, $t=\{\frac{\relspinct-p}{q}\}$ and $i(\x)+i-1=\lfloor
\frac{\relspinct-p}{q} \rfloor=\relspinc'$. The complex may be
identified with $C\{i\leq 0,k=0\}(\relspinc')$.\\

We will determine the maximum and minimum values for $\relspinct$ such
that the knot Floer homology is non-trivial, i.e.
$\widehat{\HFKT}(L,\relspinct)\neq 0$.

\begin{thm}
Suppose that $K$ is a knot in $S^3$ of genus $g(K)$, and let
$r=\frac{p}{q}\in \Q$ be a positive rational number. Under the
natural identification
$\RelSpinC(S^3,K_r)=\Z$ we will have
$$\widehat{\HFKTT}(K_r,-qg(K))\cong \widehat{\HFKTT}(K_r,qg(K)+p-1)
\cong \widehat{\HFKTT}(K,g(K))\neq 0,$$
and for any $\relspinct\in \Z$ such that $\relspinct<-qg(K)$
or $\relspinct\geq qg(K)+p$ we will have
$$\widehat{\HFKTT}(K_r,\relspinct)=0.$$
\end{thm}
\begin{proof}
Suppose that $\relspinct=-qg(K)$. Then $\Ahat_1(\relspinct)$ is
identified with the complex $C\{i\leq 0,k=0\}(\relspinc')$ where
$\relspinc'=\lfloor\frac{-qg(K)-p}{q}\rfloor<-g(K)$. This implies
that if $\x$ is a generator (intersection of $\Ta$ and $\Tb$)
such that no generator $[\x,i,0]$ is included in
$C\{i\leq 0, k=0\}(\relspinc')$ then $i(\x)<\relspinc'<-g(K)$.
Such generators will cancel each-other in homology (as
$\widehat{\HFKT}(K,\relspinc)=0$ for $\relspinc<-g(K)$, see
\cite{OS-genus}).
It is implied that the map
$$f_*:H_*(\Ahat_1(\relspinct))\lra H_*(\Bhat(\relspinct))$$
is an isomorphism. As a result we will have
$$\widehat{\HFKT}(L,-qg(K))\cong H_*(\Ahat_0(-qg(K))).$$
But $H_*(\Ahat_0(-qg(K))$ is isomorphic to $H_*(C\{i=0,k\leq
0\}(-g(K)))$. If $i(x)-k=-g(K)$ then either $i(\x)=-g(K)$ and
$k=0$, or $i(\x)<-g(K)$. The generators of the later form will
disappear in homology by the same reasoning. This implies that
$$\widehat{\HFKT}(L,-qg(K))\cong\widehat{\HFKT}(K,-g(K))\neq 0,$$
where the last non-vanishing result is borrowed from
\cite{OS-genus}.\\

It is clear that if $\relspinct<-qg(K)$ then the first isomorphism
may still be constructed. In the second part it is always implied
that $i(\x)<-g(K)$, thus the homology group
$H_*(\Ahat_0(\relspinct))$ is
trivial.\\

Now assume that $\relspinct=qg(K)+p-1$. This time
$\Ahat_0(\relspinct)$ may be identified with
$C\{i=0,k\leq 0\}(\relspinc)$ where $\relspinc=
\lfloor\frac{\relspinct}{q}\rfloor\geq g(K)$. For any
$\x\in \Ta\cap\Tb$ with $i(\x)\leq g(K)$ one can find a
non-positive integer $k$ such that $i(\x)-k=\relspinc$
(i.e. $[\x,0,k]\in C\{i=0,k\leq 0\}(\relspinc)$). As
before, this implies that
$$f_*:H_*(\Ahat_0(\relspinct))\lra H_*(\Bhat(\relspinct))=\Z$$
is an isomorphism and $\widehat{\HFKT}(L,qg(K)+p-1)\cong
H_*(\Ahat_1(qg(K)+p-1))$. Note that for this value of
$\relspinct$ we have
$$\relspinc'=\lfloor\frac{\relspinct-p}{q}\rfloor
=\lfloor\frac{qg(K)-1}{q}\rfloor=g(K)-1.$$
The generators $[\x,i,0]$ of $C\{i\leq 0,k=0\}(\relspinc')$
should then satisfy $i(\x)+i-1=g(K)-1$. This equality implies
that $i(\x)\geq g(K)$. The generators
 with $i(\x)>g(K)$ are killed in homology.  What remains is the
set of generators $[\x,0,0]$ such that $i(\x)=g(K)$ which shows that
$$\widehat{\HFKT}(L,qg(K)+p-1)\cong \widehat{\HFKT}(K,g(K))\neq 0,$$
where again we use the result of \cite{OS-genus} for the last
part. It is clear from the above argument that for $\relspinct\geq
qg(K)+p$ the knot Floer homology groups will vanish.
\end{proof}


\end{document}

%% file: command1.tex
\newcommand\Dual{\mathcal D}
\newcommand\Duality\Dual

\newcommand\RelSpinC{\underline{\SpinC}}
\newcommand\relspinc{\underline{\spinc}}

\newcommand\x{\mathbf x}

\newcommand\y{\mathbf y}

\newcommand\ModSphere{\ModFlow\left({\mathbb S}\longrightarrow
\Sym^{g-1}(\Sigma_{1})\times \Sym^2(\Sigma_{2})\right)}
\newcommand\ModSpheres\ModSphere

\newcommand\UnparModSp{\widehat \ModSp}
\newcommand\UnparModFlow\UnparModSp

\newcommand{\spinc}{\mathfrak s}

\newcommand{\spinct}{\mathfrak t}

\newcommand\ModMaps{\mathcal M}
\newcommand\ModSp\ModMaps

\newcommand\Ta{{\mathbb T}_{\alpha}}
\newcommand\Tb{{\mathbb T}_{\beta}}

\newcommand\alphas{\mbox{\boldmath$\alpha$}}

\newcommand\betas{\mbox{\boldmath$\beta$}}

\newcommand\spincrel\relspinc
\newcommand\relspinct{\underline{\mathfrak t}}

\newcommand\CFK{CFK}

\newtheorem{thm}{Theorem}[section]

\newtheorem{cor}[thm]{Corollary}

\newtheorem{lem}[thm]{Lemma}

\newtheorem{remark}[thm]{Remark}

\def\endproof{\relax\ifmmode\expandafter\endproofmath\else
  \unskip\nobreak\hfil\penalty50\hskip.75em\hbox{}\nobreak\hfil\bull
  {\parfillskip=0pt \finalhyphendemerits=0 \bigbreak}\fi}
\def\endproofmath$${\eqno\bull$$\bigbreak}
\def\bull{\vbox{\hrule\hbox{\vrule\kern3pt\vbox{\kern6pt}\kern3pt\vrule}\hrule}}

\newcommand{\Q}{\mathbb{Q}}

\newcommand{\T}{\mathbb{T}}

\newcommand{\Z}{\mathbb{Z}}

\newcommand{\ModSWfour}{\mathcal{M}}
\newcommand{\ModFlow}{\ModSWfour}

\newcommand{\SpinC}{{\mathrm{Spin}}^c}

\newcommand\abuts\Rightarrow
\newcommand\Sym{\mathrm{Sym}}


%% file: standardcommand.tex
\newcommand{\CFKT}{\text{CFK}}
\newcommand{\HFKT}{\text{HFK}}
\newcommand{\CFT}{\text{CF}}
\newcommand{\HFT}{\text{HF}}
\newcommand{\CFKTT}{\text{\emph{CFK}}}
\newcommand{\HFKTT}{\text{\emph{HFK}}}

\newcommand{\D}{\mathbb{D}}
\newcommand{\E}{\mathbb{E}}
\newcommand{\A}{\mathbb{A}}
\newcommand{\Ahat}{\widehat{\mathbb{A}}}
\newcommand{\Bhat}{\widehat{\mathbb{B}}}

\newcommand{\B}{\mathbb{B}}
\newcommand{\M}{\mathbb{M}}
\newcommand{\m}{\mathfrak{m}}

\newcommand{\gbar}{\overline{g}}
\newcommand{\fbar}{\overline{f}}
\newcommand{\xfrak}{\mathfrak{x}}
\newcommand{\yfrak}{\mathfrak{y}}
\newcommand{\Gi}{\mathcal{G}}

\newcommand{\lra}{\longrightarrow}

\newcommand{\Sig}{\Sigma}

\newcommand{\Mod}{\mathcal{M}}
